\documentclass[12pt]{amsart}
\usepackage{amscd,amsmath,amssymb,amsthm,fancyhdr}
\theoremstyle{plain}
\newtheorem{thm}{Theorem}[section]
\newtheorem{cor}[thm]{Corollary} 
\newtheorem{lemma}[thm]{Lemma} 
\newtheorem{prop}[thm]{Proposition}

\newtheorem{remark}[thm]{Remark}

\newcommand\Alg{{\operatorname{Alg}}}
\newcommand\Span{{\operatorname{Span}}}
\newcommand\card{{\operatorname{card}}}
\newcommand\diag{{\operatorname{diag}}}
\newcommand\tr{{\operatorname{tr}}}
\newcommand\Index{{\operatorname{Index}}}
\newcommand\conv{{\operatorname{conv}}}
\newcommand\id{{\operatorname{id}}}

\begin{document}

\title[Reduced Amalgamated Free Product $C^*$-Algebras]{On the Structure of Some Reduced Amalgamated Free Product $C^*$-Algebras}
\author{Nikolay A. Ivanov}
\date{\today}

\address{\hskip-\parindent
Nikolay Ivanov  \\
Department of Mathematics \\
Texas A\&M University \\
College Station, TX, 77843-3368, USA}
\email{nivanov@math.tamu.edu}

\begin{abstract}
We study some reduced free products of $C^*$-algebras with amalgamations. We give sufficient conditions for the
positive cone of the $K_0$ group to be the largest possible. We also give sufficient conditions for simplicity
and uniqueness of trace. We use the later result to give a necessary and sufficient condition for simplicity and
uniqueness of trace of the reduced $C^*$-algebras of the Baumslag-Solitar groups $BS(m,n)$. 
\end{abstract}

\maketitle

\section{Introduction}

In \cite{V85} Voiculescu indroduced the noncommutative probabilistic theory of freeness together with the 
notion of reduced amalgamated free products of $C^*$-algebras. The simplest case is amalgamation over the
complex numbers, which was considered independently by Avitzour in \cite{A82}. Avitzour also gave a sufficient condition
for simplicity and uniqueness of trace in the case of amalgamation over the complex numbers, which we
generalize here using extensively his ideas. Avitzour's work is based on the work of Powers \cite{P75}, in which Powers 
proved that the reduced $C^*$-algebra of
the free group on two generators is simple and has a unique trace. Subsequently Pashke and Salinas in
\cite{PS79} and Choi in \cite{C79} considered other reduced $C^*$-algebras of amalgams of discrete groups. 
The most general result for the case of reduced $C^*$-algebras of amalgams of discrete groups, that 
generalize Power's result is due to de la Harpe (\cite{dlH85}). It is a corollary of our result. Our result is
applicable to some HNN extensions of groups.
\par
In \cite{ABH91} Anderson, Blackadar and Haagerup studied the scale and the positive cone of ${\bf K_0}$ for the Choi algebras. In 
\cite{DR98} Dykema and R\o rdam extended their result to the case of reduced free products of
$C^*$-algebras (with amalgamation over the complex numbers). Using similar techniques we generalize the later
result for the positive cone of ${\bf K_0}$ to the case of reduced
amalgamated free products. The $K$-theory of reduced free products of nuclear $C^*$-algebras was determined by
Germain in \cite{G96} in terms of the $K$-theory of the underlying $C^*$-algebras. He gave partial results
in \cite{G97} for the $K$-theory of some reduced amalgamated free products. The question of determinig the 
$K$-theory of reduced $C^*$-algebras of amalgams of discrete groups in terms of the $K$-theory of the reduced 
$C^*$-algebras of the underlying groups was resolved completely by Pimsner in \cite{P86}.

\section{The Construction of the Reduced Amalgamated Free Product and Preliminaries} \label{S:2}

In this section we will explain the construction of reduced amalgamated free products of $C^*$-algebras of
Voiculescu, following closely \cite[\S 1]{D04}. 
\par
First we recall the definition of freeness. Suppose that we have unital $C^*$-algebras $1_{\mathfrak{A}} \in \mathfrak{B}
\subset \mathfrak{A}$ and conditional expectation $\mathfrak{E} : \mathfrak{A} \rightarrow \mathfrak{B}$.
Suppose that we have a family $\mathfrak{B} \subset \mathfrak{A}_{\iota} \subset \mathfrak{A}$, $\iota \in I$ of
$C^*$-subalgebras of $\mathfrak{A}$, all of them containing $\mathfrak{B}$. We say that the family $\{
\mathfrak{A}_{\iota} | \iota \in I \}$ is $\mathfrak{E}$-free if for any elements $a_k \in
\mathfrak{A}_{\iota_k}$, $k=1, \dots, n$, such that $\iota_1 \neq \iota_2, \iota_2 \neq \iota_3, \dots,
\iota_{n-1} \neq \iota_n$ and $\mathfrak{E}(a_k) = 0$, we have $\mathfrak{E}(a_1 a_2 \cdots a_n) = 0$. We say 
that the elements $a_i \in \mathfrak{A}, i = 1, \dots, n$ are $\mathfrak{E}$-free if the family 
$\{ C^*(\mathfrak{B} \cup \{ a_i \}) | \ i = 1, \dots, n \}$ is $\mathfrak{E}$-free. This includes the case 
$\mathfrak{B} = \mathbb{C}$ and $\mathfrak{E}$ being a state.
\par
Let $I$ be a index set, $\card(I) \geq 2$. Let $B$ be a unital $C^*$-algebra and for each $\iota \in I$ we have a
unital $C^*$-algebra $A_{\iota}$, which contains a copy of $B$ as a unital $C^*$-subalgebra. We also suppose
that for each $\iota \in I$ there is a conditional expectation $E_{\iota} : A_{\iota} \rightarrow B$,
satisfying 
\begin{equation} \label{equ:1}
\forall a \in A_{\iota},\ a \neq 0,\ \exists x \in A_{\iota},\ E_{\iota}(x^* a^* a x) \neq 0.
\end{equation}

The reduced amalgamated free product of $(A_{\iota}, E_{\iota})$ is denoted by 
$$(A, E) = \underset{\iota \in I}{*} (A_{\iota}, E_{\iota}).$$

We will be mainly interested in the case of $B \neq \mathbb{C}$ and in this case the construction depends on
some knowledge on Hilbert $C^*$-modules (see Lance's book \cite{L95} for a good exposition). 
\par
$M_{\iota}=L^2( A_{\iota}, E_{\iota})$ will denote the right Hilbert $B$-module obtained from $A_{\iota}$ by
separation and completion with respect to the norm $\| a \| = \| \langle a, a \rangle_{M_{\iota}} \|^{1/2}$,
where $\langle a_1, a_2 \rangle_{M_{\iota}} = E_{\iota}(a_1^* a_2)$. Then the linear space
$\mathcal{L}(M_{\iota})$ of all adjointable $B$-module operators on $M_{\iota}$ is actually a $C^*$-algebra
and we have a representation $\pi_{\iota} : A_{\iota} \rightarrow \mathcal{L}(M_{\iota})$ defined by
$\pi_{\iota}(a) \widehat{a'} = \widehat{aa'}$, where by $\hat{a}$ we denote the element of $M_{\iota}$,
corresponding to $a \in A_{\iota}$. $\pi_{\iota}$ is faithful by condition (\ref{equ:1}). Notice that $\pi_{\iota}|_B : B \rightarrow
\mathcal{L}(M_{\iota})$ makes $M_{\iota}$ a Hilbert $B-B$-bimodule. In this construction we have 
the specified element $\xi_{\iota} \overset{def}{=} \widehat{1_{A_{\iota}}} \in M_{\iota}$. We call the 
tripple $(\pi_{\iota}, M_{\iota}, \xi_{\iota})$ the KSGNS representation of $(A_{\iota}, E_{\iota})$, i.e. 
$(\pi_{\iota}, M_{\iota}, \xi_{\iota}) = $KSGNS$(A_{\iota}, E_{\iota})$ (KSGNS stands for Kasparov,
Steinspring, Gel$'$fand, Naimark, Segal). 
\par
For every right $B$-module $N$ one has operators $\theta_{x,y} \in \mathcal{L}(N)$ given by $\theta_{x,y}(n) =
x \langle y,n \rangle_N$ ($x,y,n \in N$). The $C^*$-subalgebra of $\mathcal{L}(N)$ that they generate is 
actually an ideal of $\mathcal{L}(N)$, which is denoted by $\mathcal{K}(N)$. It is an analogue of the
$C^*$-algebra of all compact operators on a Hilbert space. 
\par
Since for every $\iota \in I$, $\theta_{\xi_{\iota}, \xi_{\iota}} \in \mathcal{L}(M_{\iota})$ is the projection
onto the Hilbert $B-B$-subbimodule $\xi_{\iota} B$ of $M_{\iota}$ it follows that $\xi_{\iota} B$ is a
complemented submodule of $M_{\iota}$. Therefore if $P_{\iota}^{\circ} = 1 - \theta_{\xi_{\iota},
\xi_{\iota}}$ then $\pi_{\iota}(b) P_{\iota}^{\circ} = P_{\iota}^{\circ} \pi_{\iota}(b) \in
\mathcal{L}(M_{\iota})$ for each $b \in B$. We define $M_{\iota}^{\circ} \overset{def}{=} P_{\iota}^{\circ} M_{\iota}$. If we view $\xi
\overset{def}{=} 1_B$ as an element of the Hilbert $B-B$-bimodule $B$, we can define
\begin{equation} \label{equ:2}
M = \xi B \oplus \underset{\iota_1 \neq \iota_2, \iota_2 \neq \iota_3, \dots, \iota_{n-1} \neq \iota_n}{\underset{\iota_1, \dots, \iota_n
\in I}{\underset{n\in \mathbb{N}}{\bigoplus}}} M_{\iota_1}^{\circ} \otimes_B M_{\iota_2}^{\circ} \otimes_B \cdots \otimes_B 
M_{\iota_n}^{\circ},
\end{equation}
where $\otimes_B$ means interior tensor product (see \cite{L95}). The Hilbert $B-B$-bimodule $M$ constructed above is
called the free product of $\{ M_{\iota}, \iota \in I \}$ with respect to vectors $\{ \xi_{\iota}, \iota \in I \}$ and is denoted by $(M,
\xi) = \underset{\iota \in I}{*}(M_{\iota}, \xi_{\iota}).$ 
\par
For each $\iota \in I$ set 
\begin{equation} \label{equ:3}
M(\iota) = \eta_{\iota} B \oplus \underset{\iota_1 \neq \iota}{\underset{\iota_1 \neq \iota_2, \iota_2 \neq \iota_3, \dots, \iota_{n-1} \neq \iota_n}{\underset{\iota_1, \dots, \iota_n
\in I}{\underset{n\in \mathbb{N}}{\bigoplus}}}} M_{\iota_1}^{\circ} \otimes_B M_{\iota_2}^{\circ} \otimes_B \cdots \otimes_B 
M_{\iota_n}^{\circ},
\end{equation}
where $\eta_{\iota} \overset{def}{=} 1_B \in B$. We define a unitary operator $$V_{\iota} : M_{\iota} \otimes_B M(\iota) \rightarrow M$$
given on elementary tensors by:
\begin{align*}
[& \xi_{\iota}] \otimes [\eta_{\iota}] \mapsto \xi,  \\ 
[& \zeta] \otimes [\eta_{\iota}] \mapsto \zeta, \text{ where $\zeta \in M^{\circ}_{\iota} \subset M$}  \\
[& \xi_{\iota}] \otimes [\zeta_1 \otimes \cdots \otimes \zeta_n] \mapsto \zeta_1 \otimes \cdots \otimes 
\zeta_n, \text{ where $\zeta_j \in M^{\circ}_{\iota_j}$ and } \\ 
  &  \hspace{3 in}  \iota \neq \iota_1, \iota_1 \neq \iota_2, \dots, \iota_{n-1} \neq  \iota_n   \\ 
[& \zeta] \otimes [\zeta_1 \otimes \cdots \otimes \zeta_n] \mapsto \zeta \otimes \zeta_1 \otimes \cdots \otimes \zeta_n, \text{ where $\zeta
\in M^{\circ}_{\iota}$ and}  \\
 &  \hspace{2.5 in} \zeta_j \in M^{\circ}_{\iota_j} \text{ with } \iota \neq \iota_1, \iota_1 \neq \iota_2, \dots, \iota_{n-1} \neq \iota_n.
\end{align*}
Let $\lambda_{\iota} : A_{\iota} \rightarrow \mathcal{L}(M)$ be the $*$-homomorphism given by
$\lambda_{\iota}(a) = V_{\iota}(\pi_{\iota}(a)
\otimes 1)V_{\iota}^*$. $\lambda_{\iota}$ defines a left action of $A_{\iota}$ on $M$. Condition (\ref{equ:1}) implies that $\lambda_{\iota}$ is injective. Then $A$ is defined as the $C^*$-subalgebra of $\mathcal{L}(M)$, generated by $\underset{\iota \in I}{\cup}
\lambda_{\iota}(A_{\iota})$, and $E: A \rightarrow B$ is the conditional expectation, given by $E(a) = \langle \xi, a(\xi) \rangle_M$.
Note that if $b \in B$, then $\lambda_{\iota}(b) \in \mathcal{L}(M)$ does not depend on $\iota$.
$\lambda_{\iota}(b)$ gives the left action of $B$ on $M$. Because of
condition (\ref{equ:1}) for each $\iota \in I$ we have unital embeddings $A_{\iota} \hookrightarrow A$, which come from the $*$-homorphisms $\lambda_{\iota} :
A_{\iota} \rightarrow \mathcal{L}(M)$. We will denote by $\pi$ the representation $\pi : A \rightarrow
\mathcal{L}(M)$ arising from the reduced amalgamated free product construction. We actually have that $(\pi,
M, \xi) = \text{KSGNS}(A, E)$. 
\par
Set $A_{\iota}^{\circ} = A_{\iota} \cap \ker(E_{\iota})$. For $a \in A_{\iota}^{\circ},\ \zeta_j \in M^{\circ}_{\iota_j}$ with $\iota_1,
\dots, \iota_n \in I, n \geq 2$, and $\iota_j \neq \iota_{j+1}$ we have 
\begin{equation} \label{equ:4}
\lambda_{\iota}(a)(\zeta_1 \otimes \cdots \otimes \zeta_n) = \begin{cases}
                                                             \widehat{a} \otimes \zeta_1 \otimes \cdots \otimes \zeta_n, & \text{ 
                                                             if } \iota \neq \iota_1, \\
                                                             (a(\zeta_1) - \xi_{\iota_1} \langle
                                                             \xi_{\iota_1}, a(\zeta_1) \rangle) \otimes
                                                             \zeta_2 \otimes \cdots \otimes \zeta_n + & \\ 
                                                             \hspace{1in} \pi_{\iota_2}(\langle \xi_{\iota_1},
                                                             a(\zeta_1) \rangle) \zeta_2 \otimes \cdots \otimes \zeta_n, & \text{  if } 
                                                             \iota = \iota_1.
                                                             \end{cases}
\end{equation}

We will omit writing $\lambda_{\iota}$ and $\pi_{\iota}$ if this leads to no confusion.
\par
We will use the following notation for the case of amalgamation which is similar to the notation in \cite{D99LN}
used for the case of amalgamation over the scalars. If everything is as above by $\Lambda_B^{\circ}(\{ A_{\iota}^{\circ}
| \iota \in I \})$ we will denote the set of words of the form $a_1 a_2 \cdots a_n$, where $n \geq 1$ and 
$a_j \in A_{\iota_j}^{\circ}$ with $\iota_j \neq \iota_{j+1}$ for $1 \leq j \leq n-1$. We will not distinguish
between two words from $\Lambda_B^{\circ}(\{ A_{\iota}^{\circ} | \iota \in I \})$ which are equal as elements
of $A$. We will denote $\Lambda_B(\{ A_{\iota}^{\circ} | \iota \in I \}) \overset{def}{=} B \cup \Lambda_B^{\circ}(\{ A_{\iota}^{\circ}
| \iota \in I \})$. By $\mathbb{C}(A)$ we will denote the span of words from $\Lambda_B(\{ A_{\iota}^{\circ} |
\iota \in I \})$. Notice that $\mathbb{C}(A)$ is norm-dense in $A$. For a word $a_1 a_2 \cdots a_n \in
\Lambda_B^{\circ}(\{ A_{\iota}^{\circ} | \iota \in I \})$, where $ n \geq 1$, $a_j \in A_{\iota_j}^{\circ}$ 
with $\iota_j \neq \iota_{j+1}$ for $1 \leq j \leq n-1$ we will consider to be of length $n$. Elements of $B$ we will
consider to be of length $0$. 
\par
We will be mainly interested in the case $\card(I) = 2$ and that there exist states $\phi_{\iota}$ on $A_{\iota}$ for $\iota = 1,2$,
such that these states are invariant under $E_{\iota}$, i.e. for $\iota = 1,2$ and $\forall a_{\iota} \in A_{\iota}$ we have
$\phi_{\iota}(a_{\iota}) = \phi_{\iota}(E_{\iota}(a_{\iota}))$. We also require $\phi_1(b) = \phi_2(b)$ for $b \in B$. $\phi
\overset{def}{=} \phi_B \circ E$, where $\phi_B  \overset{def}{=} \phi_1|_B = \phi_2|_B$ is a well defined $E$-invariant state on
$(A,E) = (A_1,E_1) * (A_2, E_2)$. In such case we will write formally $$(A, E, \phi) = (A_1, E_1, \phi_1) * 
(A_2, E_2, \phi_2),$$ although the construction of $(A, E)$ does not depend on $\phi_{\iota}, \iota =1,2$.  
\par
Let $Q_{\iota} : M \to M_{\iota}$ be the orthogonal projection of $M$ onto the the complemented submodule
$M_{\iota} = M_{\iota}^{\circ} \oplus B \widehat{1_A}$ (see (\ref{equ:2}). It is easy to see that $F_{\iota} : A
\to A_{\iota}$ defined by $F_{\iota}(a) = \lambda_{\iota}(Q_{\iota} a Q_{\iota}^*)$ for $a \in A$ is a 
conditional expectation from $A$ onto $A_{\iota} \subset A$ which is invariant with respect to $E$, i.e. $E(a) =
E(F_{\iota}(a))$ for all $a \in A$.  
\par                                 
We will need the following result concerning the faithfulness of $E$. This short proof was noted to me by
$\acute{E}$ric Ricard:

\begin{thm}[$\acute{E}$. Ricard] \label{thm:faithfulness}
Let $(A, E) = \underset{\iota \in I}{*} (A_{\iota}, E_{\iota})$. Then the faithfulness of $E_{\iota}$ for
$\forall \iota \in I$ implies the faithfulness of $E$.
\end{thm}

\begin{proof}
Let $K=\ker(E) \cap A^+$, then for any $x \in K$, $F_{\iota}(x)=0$ for each $\iota \in I$ because of 
$0=E(x)=E(F_{\iota}(x))$ and the faithfulness of $E_{\iota} = E|_{A_{\iota}}$. 
Then $x \in K$ implies $a_{\iota}^*xa_{\iota} \in K$ for each $a_{\iota} \in A_{\iota}$, since 
$E(a_{\iota}^*xa_{\iota})=E(F_{\iota}(a_{\iota}^*xa_{\iota}))=E(a_{\iota}^*F_{\iota}(x)a_{\iota})=0$.
\par
Now if $x \in A$ is such that $x^*x \in K$, then for $y=x  a_1 \cdots a_n$ we have $y^*y\in K$, where 
$a_1 a_2 \cdots a_n \in \Lambda_B^{\circ}(\{ A_{\iota}^{\circ} | \iota \in I \})$. 
In particular, it means that
$$0 = E(y^*y) = E((x  a_1 \cdots a_n)^* x  a_1 \cdots a_n) = \langle \widehat{1_A}, (x  a_1
\cdots a_n)^*x  a_1 \cdots a_n(\widehat{1_A}) \rangle_M =$$ 
$$= \langle x a_1 \cdots a_n(\widehat{1_A}), 
(x a_1 \cdots a_n)^* x a_1 \cdots a_n(\widehat{1_A}) \rangle_M = \langle x(\widehat{a_1} \otimes \cdots
\otimes \widehat{a_n}) , x(\widehat{a_1} \otimes \cdots
\otimes \widehat{a_n}) \rangle_M$$ 
so $x$ vanishes on the dense subset $\widehat{\mathbb{C}(A)}$ of $M$. 
Therefore $x \equiv 0$ as an operator of $\mathcal{L}(M)$. Thus $K = \{ 0 \}$.
\par
This proves the theorem.
\end{proof}

It follows immediately from this theorem that if $(A, E, \phi) = (A_1, E_1, \phi_1) * (A_2, E_2, \phi_2)$ then the faithfulness of
$\phi_1$ and $\phi_2$ imply the faithfulness of $\phi$.
\par
Now let's define the sets
\begin{equation} \label{equ:5}
\Lambda_B^1 \overset{def}{=} \ \Span(\underset{k=0}{\overset{\infty}{\bigcup}} A_1^{\circ} (A_2^{\circ}
A_1^{\circ})^k ) \subset \mathbb{C}(A)
\end{equation}
and
\begin{equation} \label{equ:6}
\Lambda_B^2 \overset{def}{=} \ \Span(\underset{k=0}{\overset{\infty}{\bigcup}} A_2^{\circ} (A_1^{\circ}
A_2^{\circ})^k ) \subset \mathbb{C}(A).
\end{equation}

Define also
\begin{equation*}
\Lambda_B^{21} \overset{def}{=} \ \Span(\underset{k=1}{\overset{\infty}{\bigcup}} (A_2^{\circ}
A_1^{\circ})^k ) \subset \mathbb{C}(A)
\end{equation*}
and
\begin{equation*}
\Lambda_B^{12} \overset{def}{=} \ \Span(\underset{k=1}{\overset{\infty}{\bigcup}} (A_1^{\circ}
A_2^{\circ})^k ) \subset \mathbb{C}(A).
\end{equation*}

Some of the most important examples are those of reduced $C^*$-algebras of amalgams of discrete groups. For each discrete group $N$ we have the canonical
tracial state $\tau_N \overset{def}{=} \
\langle \cdot, \widehat{1_H} \rangle_{l^2(H)}$ on $C^*_r(N)$. For
each subgroup $S$ of $N$ we have a canonical conditional expectation $E_S^N : C^*_r(N) \rightarrow C^*_r(S)$ given on elements 
$\{ \lambda_{n}, n \in N \}$ by 
$$E_S^N(\lambda_n) = \begin{cases}
                     \lambda_n, & \text{  if } n \in S, \\
                     0, & \text{  if } n \notin S.
                     \end{cases}$$

Let $G_1 \supset H \subset G_2$ be two discrete groups, containing a common subgroup (an isomorphic copy of $H$). Then we have 
$(C^*_r(G), E_H^{G}) =  (C^*_r(G_1), E_H^{G_1}) * (C^*_r(G_2), E_H^{G_2})$, where $G = G_1 \underset{H}{*} G_2$. 
\par
The canonical tracial states $\tau_{G_{\iota}}, \iota =1,2$ and $\tau_G$ are invariant under $E_H^{G_{\iota}}, \iota =1,2$ and $E_H^G$
respectivelly and $\tau_G = \tau_H \circ E_H^G$. Thus we can write formally 
$$(C^*_r(G), E_H^{G}, \tau_G) = (C^*_r(G_1), E_H^{G_1}, \tau_{G_1}) * (C^*_r(G_2), E_H^{G_2}, \tau_{G_2}).$$ 

\section{${\bf K_0}^+$}

We give the results of Germain and Pimsner first.

\begin{thm}[\cite{G96}] \label{thm:germain}
Let $(A, \phi) = (A_1, \phi_1) * (A_2, \phi_2)$ is the reduced free product (with amalgamation over
$\mathbb{C}$) of the unital, nuclear $C^*$-algebras $A_1$ and $A_2$ with respect to states $\phi_1$ and 
$\phi_2$. Then we have the following six term exact sequence:
\begin{equation*} 
 \begin{CD}
  \mathbb{Z} \cong & {\bf K_0}(\mathbb{C}) @>({\bf K_0}(i_1), -{\bf K_0}(i_2))>> {\bf K_0}(A_1) \oplus {\bf K_0}(A_2) 
    @>{\bf K_0}(j_1)+{\bf K_0}(j_2)>> {\bf K_0}(A)  & \\
              &    @AAA            & &                             @VVV        & \\
   & {\bf K_1}(A) @<{\bf K_1}(j_1)+{\bf K_1}(j_2)<< {\bf K_1}(A_1) \oplus {\bf K_1}(A_2) 
    @<({\bf K_1}(i_1), -{\bf K_1}(i_2))<< {\bf K_1}(\mathbb{C}) & \cong 0, \\
 \end{CD}
\end{equation*}
where $i_k : \mathbb{C} \rightarrow A_k$ are the the unital $*$-homorphisms and $j_k : A_k \rightarrow A$ are
the unital embeddings arising from the construction of reduced free product ($k=1,2$).
\end{thm}

\begin{thm}[\cite{P86}] \label{thm:pimsner}
Suppose that $G_1 \supset H \subset G_2$ are countable, discrete groups. Let $G=G_1 \underset{H}{*} G_2$. Then we have the following six term
exact sequence:
\begin{equation*} 
 \begin{CD}
  {\bf K_0}(C^*_r(H)) @>({\bf K_0}(i_1), -{\bf K_0}(i_2))>> {\bf K_0}(C^*_r(G_1)) \oplus 
  {\bf K_0}(C^*_r(G_2)) @>{\bf K_0}(j_1)+{\bf K_0}(j_2)>> {\bf K_0}(C^*_r(G))  \\
                 @AAA            & &                             @VVV        \\
  {\bf K_1}(C^*_r(G)) @<{\bf K_1}(j_1)+{\bf K_1}(j_2)<< {\bf K_1}(C^*_r(G_1)) \oplus {\bf K_1}(C^*_r(G_2)) 
    @<({\bf K_1}(i_1), -{\bf K_1}(i_2))<< {\bf K_1}(C^*_r(H)), \\
 \end{CD}
\end{equation*}
where $i_k : C^*_r(H) \rightarrow C^*_r(G_k)$ and $j_k : C^*_r(G_k) \rightarrow C^*_r(G)$ are the canonical
inclusion maps ($k=1,2$).
\end{thm}

Now suppose that we have unital $C^*$-algebaras $A_{\iota},\ \iota=1,2$ and $B$. Suppose that we have unital
inclusions $B \hookrightarrow A_{\iota}$ and conditional expectations $E_{\iota} : A_{\iota} 
\rightarrow B$ that satisfy property (\ref{equ:1}). Suppose also that for $\iota = 1,2$ we have tracial states
$\tau_{\iota}$ on $A_{\iota}$ which satisfy $\tau_B \overset{def}{=} \tau_1|_B = \tau_2|_B$ and which
are invariant under $E_{\iota}$, i.e $\tau_{\iota}(a_{\iota}) = \tau_{\iota}(E_{\iota}(a_{\iota}))$ for each
$a_{\iota} \in A_{\iota}$. Let us denote $(A, E, \tau) \overset{def}{=} (A_1, E_1, \tau_1) * (A_2, E_2,
\tau_2)$ and let $j_{\iota} : A_{\iota} \rightarrow A$ are the inclusion maps, coming from the construction of
reduced amalgamated free products. Suppose that $\tau \overset{def}{=} \tau_B \circ E$ is a faithful tracial state. Let's define 
$$\Gamma \overset{def}{=} {\bf K_0}(j_1)({\bf K_0}(A_1)) + {\bf K_0}(j_2)({\bf K_0}(A_2))
\subset {\bf K_0}(A).$$ 

Then every element in $\Gamma$ can be represented as 
\begin{equation*} 
([p_1]_{{\bf K_0}(A)} - [q_1]_{{\bf K_0}(A)}) + ([p_2]_{{\bf K_0}(A)} - [q_2]_{{\bf K_0}(A)}),
\end{equation*}
where $p_{\iota}, q_{\iota}$ are projections in some matrix algebras over $A_{\iota}$ for $\iota = 1,2$. By
expanding those matrices and adding zeros we can suppose without loss of generality that $p_{\iota},
q_{\iota}$ are projections from $M_n(A_{\iota})$ for some $n \in \mathbb{N}$ for $\iota = 1,2$. Therefore
every element of $\Gamma$ can be represented in the form
\begin{equation} \label{equ:7}
\left[ \left( \begin{array}{cc} p_1 & 0 \\ 0 & p_2 \end{array} \right) \right]_{{\bf K_0}(A)} - 
\left[ \left( \begin{array}{cc} q_1 & 0 \\ 0 & q_2 \end{array} \right) \right]_{{\bf K_0}(A)}, 
\end{equation}
where now 
\begin{equation*} 
\left( \begin{array}{cc} p_1 & 0 \\ 0 & p_2 \end{array} \right) \text{ and } 
\left( \begin{array}{cc} q_1 & 0 \\ 0 & q_2 \end{array} \right) \ \in M_{2n}(A). 
\end{equation*}

We want to obtain a sufficient condition so that all elements $\gamma \in \Gamma$ for which ${\bf K_0}(\tau)(\gamma) > 0$
come from projections, i.e. $\exists m \in \mathbb{N}$ and $\rho \in M_m(A)$, such that 
$\gamma = [ \rho ]_{{\bf K_0}(A)}$ in ${\bf K_0}(A)$. 
\par
By definition the positive cone of ${\bf K_0}(A)$ is
$${\bf K_0}(A)^+ = \{ x \in {\bf K_0}(A) | \exists p \text{ projection in $M_n(A)$ for some $n$ with } 
x = [p]_{{\bf K_0}(A)} \}.$$

The scale of ${\bf K_0}(A)$ is
$$\Sigma(A) = \{ x \in {\bf K_0}(A) | \exists p \text{ projection in $A$ with } 
x = [p]_{{\bf K_0}(A)} \}.$$ 

Dykema and R\o rdam proved the following:

\begin{thm}[\cite{DR98}] \label{thm:dr}
Let $(A, \tau) = (A_1, \tau_1) * (A_2, \tau_2)$ be the reduced free product of the unital $C^*$-algebras $A_1$
and $A_2$ with respect to the faithful tracial states $\tau_1$ and $\tau_2$. Suppose that the Avitzour
condition holds, namely there exist unitaries $u_1 \in A_1$ and $u_2, u_2' \in A_2$, such that $\tau_1(u_1) =
\tau_2(u_2) = \tau_2(u_2') = \tau_2(u_1^* u_1') = 0$. Then we have 
$$\Gamma \cap {\bf K_0}(A)^+ = \{ \gamma \in \Gamma | {\bf K_0}(\tau)(\gamma) > 0 \} \cup \{ 0 \}$$
and
$$\Gamma \cap \Sigma(A) = \{ \gamma \in \Gamma | 0 < {\bf K_0}(\tau)(\gamma) < 1 \} \cup \{ 0,1 \}.$$
\end{thm}

Notice that Theorem \ref{thm:germain} implies that if $A_1$ and $A_2$ are nuclear then $\Gamma = {\bf
K_0}(A)$.
\par
Anderson, Blackadar and Haagerup proved this theorem for the case of $A = C^*_r(\mathbb{Z}_n * \mathbb{Z}_m)$
and gave one of the main technical tool for proving Theorem \ref{thm:dr}, which we will use here also:

\begin{prop}[\cite{ABH91}] \label{prop:abh}
Let $\mathfrak{A}$ be a unital $C^*$-algebra and let $\phi$ be a faithful state on  $\mathfrak{A}$. Suppose
that $p,q \in \mathfrak{A}$ are projections that are $\phi$-free in $\mathfrak{A}$. If $\phi(p) < \phi(q)$
then $\| p(1-q) \| < 1$ and there is a partial isometry $\nu \in \mathfrak{A}$ such that $\nu \nu^* = p$ and
$\nu^* \nu < q$.
\end{prop}

Now we can state and prove our result:

\begin{thm} \label{thm:crit}
Let $A_{\iota}$ be unital $C^*$-algebras that contain the unital $C^*$-algebra $B$ as a unital
$C^*$-subalgebra, i.e. $1_{A_{\iota}} \in B \subset A_{\iota}$, $\iota=1,2$. Suppose that we have conditional expectations
$E_{\iota} : A_{\iota} \rightarrow B$ and faithful tracial states $\tau_{\iota}$ on $A_{\iota}$ for $\iota
=1,2$ such that $\tau_{\iota} = \tau_{\iota} \circ E_{\iota}$ and $\tau_1|_B = \tau_2|_B$. Form
the reduced amalgamated free product $(A, E, \tau) = (A_1, E_1, \tau_1) * (A_2, E_2, \tau_2)$. Suppose that
the following two conditions hold:
\begin{equation} \label{equ:8}
         \begin{cases}
         \forall b_1, \dots, b_l \in B, \text{ with } \tau(b_1) = \dots = \tau(b_l) =0,\ \exists m \in
         \mathbb{N} \text{ and unitaries }  \\
         \nu_{11}, \dots, \nu_{1m}, \nu_{21}, \dots, \nu_{2m} \text{ such that } \nu_{12}, \dots, \nu_{1m} \in A_1^{\circ},\ \nu_{21}, \dots,
         \nu_{2(m-1)} \in A_2^{\circ}, \text{ and: } \\
         \text{ either } \nu_{11} \in A_1^{\circ},\ \nu_{2m} \in A_2^{\circ} \text{ or } \\ 
         \nu_{11} = 1_{A_1},\ \nu_{2m} \in A_2^{\circ}, \text{ or } \\
         \nu_{11} \in A_1^{\circ},\ \nu_{2m} = 1_{A_2}, \\ 
         \nu_{11} = 1_{A_1},\ \nu_{2m} = 1_{A_2},\ k \geq 2 \\ 
         \text{ with } E((\nu_{11} \nu_{21} \nu_{12} \cdots \nu_{1m} \nu_{2m}) b_k (\nu_{11} \nu_{21} \nu_{12}
         \cdots \nu_{1m} \nu_{2m})^*) = 0 \text{ for } k=1, \dots, l, \\ 
         (\text{i.e. there are unitaries that conjugate } B \ominus \mathbb{C} 1_B \text{ out of } B) \\ 
         \end{cases}
\end{equation}
and 
\begin{equation} \label{equ:9}
\exists \text{ unitaries } u_1 \in A_1^{\circ},\ u_2, u_2' \in A_2^{\circ}, \text{ with } E_2(u_2 u_2'^*)=0.
\end{equation}

Then:
\begin{equation} \label{equ:10}
\Gamma \cap {\bf K_0}(A)^+ = \{ \gamma \in \Gamma | {\bf K_0}(\tau)(\gamma) > 0 \} \cup \{ 0 \}.
\end{equation}
\end{thm}

\begin{proof}
All elements of $\Gamma$ have the form (\ref{equ:7}) for some  $n \in \mathbb{N}$ and projections $p_1, q_1$ from 
$M_n(A_1)$ and $p_2, q_2$ from $M_n(A_2)$. Denote
\begin{equation*} 
\gamma = \left[ \left( \begin{array}{cc} p_1 & 0 \\ 0 & p_2 \end{array} \right) \right]_{{\bf K_0}(A)} - 
\left[ \left( \begin{array}{cc} q_1 & 0 \\ 0 & q_2 \end{array} \right) \right]_{{\bf K_0}(A)}. 
\end{equation*}

Consider
$$P \overset{def}{=} \left( \begin{array}{cc} U_2 & 0 \\ 0 & U_2 U_1 \end{array} \right)
\left( \begin{array}{cc} p_1 & 0 \\ 0 & p_2 \end{array} \right) 
\left( \begin{array}{cc} U_2^* & 0 \\ 0 & U_1^* U_2^* \end{array} \right) \text{ and }$$  
$$Q \overset{def}{=} \left( \begin{array}{cc} U_2 & 0 \\ 0 & U_2 U_1 \end{array} \right)
\left( \begin{array}{cc} q_1 & 0 \\ 0 & q_2 \end{array} \right)  
\left( \begin{array}{cc} U_2^*  & 0 \\ 0 & U_1^* U_2^*  \end{array} \right),$$
where $U_1 = \diag(u_1, \dots, u_1) \in M_n(A_1)$ and $U_2 = \diag(u_2, \dots, u_2) \in M_n(A_2)$. \\
It is clear that $P, Q \in M_{2n}(\Lambda_B^2 \oplus B 1_B)$. 
For $T \in M_m(A)$ we will denote by $T_{ij}$ the $ij$-entry of
$T$. Now consider the set of elements $S_P = \{ E(P_{ij}) - \tau(P_{ij}) | 1 \leq i,j \leq 2n \} \cup \{ E(u_1
P_{ij} u_1^*) - \tau(u_1 P_{ij} u_1^*) | 1 \leq i,j \leq 2n \}$ and the set
$S_Q = \{ E(Q_{ij}) - \tau(Q_{ij}) | 1 \leq i,j \leq 2n \} \cup \{ E(u_1 Q_{ij} u_1^*) - \tau(u_1 Q_{ij} u_1^*) | 1 
\leq i,j \leq 2n \}$. 
\par
Applying condition (\ref{equ:8}) to the set $S_P$ we obtain unitaries $\nu_{ij}, i=1,2, j=1, \dots,
m_P$. 
\par
Set 
$$W_P \overset{def}{=} \begin{cases}
                       \nu_{11} \nu_{21} \nu_{12} \cdots \nu_{2 (m_{P}-1)} \nu_{1m_P},& \text{  if } \nu_{2m_P} = 1_{A_2},\ 
                       \nu_{11} \in A_1^{\circ}, \\ 
                       \nu_{11} \nu_{21} \nu_{12} \cdots \nu_{2 (m_P - 1)} \nu_{1m_P} \nu_{2m_P} u_1,& \text{  if } \nu_{2m_P}
                       \in A_2^{\circ} ,\ \nu_{11} \in A_1^{\circ}, \\ 
                       u_1 \nu_{21} \nu_{12} \cdots \nu_{2 (m_P - 1)} \nu_{1m_P} \nu_{2m_P} u_1,& \text{  if } \nu_{2m_P} \in A_2^{\circ},\ 
                       \nu_{11} = 1_{A_1}. \\
                       u_1 \nu_{21} \nu_{12} \cdots \nu_{2 (m_P - 1)} \nu_{1m_P},& \text{  if } \nu_{2m_P} =
                       1_{A_2},\ \nu_{11} = 1_{A_1},\ k \geq 2. \\
                       \end{cases}$$
                       
Applying condition (\ref{equ:8}) to the set $S_Q$ we obtain unitaries $\nu_{ij}', i=1,2, j=1, \dots, m_Q$.
\par 
Set
$$W_Q \overset{def}{=} \begin{cases}
                       \nu'_{11} \nu'_{21} \nu'_{12} \cdots \nu'_{2 (m_{P}-1)} \nu'_{1m_P},& \text{  if } \nu'_{2m_P} = 1_{A_2},\ 
                       \nu'_{11} \in A_1^{\circ}, \\ 
                       \nu'_{11} \nu'_{21} \nu'_{12} \cdots \nu'_{2 (m_P - 1)} \nu'_{1m_P} \nu'_{2m_P} u_1,&
                       \text{  if } \nu'_{2m_P}
                       \in A_2^{\circ} ,\ \nu'_{11} \in A_1^{\circ}, \\ 
                       u_1 \nu'_{21} \nu'_{12} \cdots \nu'_{2 (m_P - 1)} \nu'_{1m_P} \nu'_{2m_P} u_1,& \text{ 
                       if } \nu'_{2m_P} \in A_2^{\circ},\ 
                       \nu'_{11} = 1_{A_1}. \\
                       u_1 \nu'_{21} \nu'_{12} \cdots \nu'_{2 (m_P - 1)} \nu'_{1m_P},& \text{  if } \nu'_{2m_P} =
                       1_{A_2},\ \nu'_{11} = 1_{A_1},\ k \geq 2. \\
                       \end{cases}$$

It is easy to see that $W_P P W_P^*,\ W_Q Q W_Q^* \in M_{2n}(\Lambda_B^1 \oplus \mathbb{C} 1_B)$.
\par
Now consider the following matrix in $M_{2n}(A)$: 
$$U = (\frac{\omega^{ij}}{\sqrt{2n}} u_2' (u_1 u_2)^{2ni + j} u_2'^*)_{i,j = 1}^{2n},$$
where $\omega = \exp(2 \pi \sqrt{-1} / 2n)$ is a primitive $2n$-th root of $1$. 
It is clear that $U \in M_{2n}(\Lambda_B^2)$. We will check that $U$ is a unitary matrix: 
$$(UU^*)_{ij} = (2n)^{-1} \underset{k=1}{\overset{2n}{\sum}}
\omega^{ik} u_2' (u_1 u_2)^{2ni+k} \omega^{-jk} (u_1 u_2)^{-2nj-k} u_2'^* =$$
$$(2n)^{-1} \underset{k=1}{\overset{2n}{\sum}} \omega^{(i-j)k} u_2' (u_1 u_2)^{2n(i-j)} u_2'^* = (2n)^{-1} u_2' (u_1
u_2)^{2n(i-j)} u_2'^* \underset{k=1}{\overset{2n}{\sum}} \omega^{(i-j)k} = \delta_{ij} 1_A.$$
$$(U^*U)_{ij} = (2n)^{-1} \underset{k=1}{\overset{2n}{\sum}}
\omega^{-ik} u_2' (u_1 u_2)^{-2nk-i} \omega^{jk} (u_1 u_2)^{2nk+j} u_2'^* =$$ 
$$(2n)^{-1} \underset{k=1}{\overset{2n}{\sum}} \omega^{(j-i)k} u_2' (u_1 u_2)^{j-i} u_2'^* = (2n)^{-1} u_2' (u_1
u_2)^{j-i} u_2'^* \underset{k=1}{\overset{2n}{\sum}} \omega^{(j-i)k} = \delta_{ij} 1_A.$$

Thus $U \in M_{2n}(A)$ is a unitary. 
\par
Take $T \in M_{2n}(\Lambda_B^1 \oplus \mathbb{C} 1_B)$. Then $T = T_0 + T_1 \otimes 1_A$, with $T_0 \in M_{2n}(\Lambda_B^1)$ and $T_1
\in M_{2n}(\mathbb{C})$. It is easy to see that $U T_0 U^* \in M_{2n}(\Lambda_B^2)$. Now if $T_1 =
(t_{ij})_{i,j=1}^{2n}$ then for $U (T_1 \otimes 1_A) U^* = (s_{ij})_{i,j=1}^{2n}$ we have 
$$s_{ij} = (2n)^{-1} \underset{k=1}{\overset{2n}{\sum}} \underset{l=1}{\overset{2n}{\sum}} \omega^{ik} u_2'
(u_1 u_2)^{2ni+k} u_2'^* t_{kl} \omega^{-jl} u_2' (u_1 u_2)^{-2nj-l} u_2'^* = $$ 
$$(2n)^{-1} \underset{k=1}{\overset{2n}{\sum}} \underset{l=1}{\overset{2n}{\sum}} t_{kl} \omega^{ik-jl} u_2' (u_1
u_2)^{2ni+k-2nj-l} u_2'^*.$$ 

If $i \neq j$ then $2ni+k-2nj-l \neq 0$ for any $1 \leq k,l \leq 2n$, so in this case $s_{ij} \in \Lambda_B^2$. If
$i=j$ then:
$$s_{ii} = (2n)^{-1} \underset{k=1}{\overset{2n}{\sum}} \underset{l=1}{\overset{2n}{\sum}} t_{kl}
\omega^{i(k-l)} u_2' (u_1 u_2)^{k-l} u_2'^* = $$
$$(2n)^{-1} \underset{k \neq l}{\underset{1\leq k,l \leq 2n}{\sum}} 
t_{kl} \omega^{i(k-l)} u_2' (u_1 u_2)^{k-l} u_2'^* + ((2n)^{-1} \underset{k=1}{\overset{2n}{\sum}} t_{kk})
\otimes 1_A.$$

So $s_{ii} = s'_{ii} + \tr_{2n}(T_1) \otimes 1_A$, where $s'_{ii} \in \Lambda_B^2$. All this means that $U (T_1 \otimes 1_A) U^*
= T_1' + \tr_{2n}(T_1) 1_A \otimes 1_{M_{2n}(\mathbb{C})}$, with $T_1' \in M_{2n}(\Lambda_B^2)$, which implies that 
$U T U^* \in M_{2n}(\Lambda_B^2) \oplus \mathbb{C} 1_{M_{2n}(A)}$. 
\par
This means that we have
\begin{equation} \label{equ:11}
P' \overset{def}{=} U W_P P W_P^* U^* \in M_{2n}(\Lambda_B^2) \oplus \mathbb{C} 1_{M_{2n}(A)}
\end{equation}
and 
\begin{equation} \label{equ:12}
Q' \overset{def}{=} u_1 U W_Q Q W_Q^* U^* u_1^* \in M_{2n}(\Lambda_B^1) \oplus \mathbb{C} 1_{M_{2n}(A)}.
\end{equation}

It is clear that $\tr_{2n} \otimes E (P') = \tr_{2n} \otimes \tau (P')$ and
that $\tr_{2n} \otimes E (Q') = \tr_{2n} \otimes \tau (Q')$. Since $P'$ and $Q'$ are nontrivial projections 
it is also clear that $C^*(\{ P', 1_A \})$ and $C^*(\{ Q', 1_A \})$ are both $2$-dimensional. 
Therefore for any $p \in C^*(\{ P', 1_A \})$ and
$q \in C^*(\{ Q', 1_A \})$ we have $\tr_{2n} \otimes E (p) = \tr_{2n} \otimes \tau (p)$ and $\tr_{2n} \otimes
E (q) = \tr_{2n} \otimes \tau (q)$. Therefore from (\ref{equ:11}), (\ref{equ:12}) and the definition of freeness it follows that 
$P'$ is both $\tr_{2n} \otimes E$-free and $\tr_{2n} \otimes \tau$-free from $Q'$. 
\par
Since $\tr_{2n} \otimes \tau$ is a faithful tracial state (because of faithfulness of $\tau_1$, $\tau_2$ and
Theorem \ref{thm:faithfulness}) and because
$$\tr_{2n} \otimes \tau(P') = (2n)^{-1} {\bf K_0}(\tau)(P) > (2n)^{-1} {\bf
K_0}(\tau)(Q) = \tr_{2n} \otimes \tau(Q'),$$ 
we can apply Proposition \ref{prop:abh} and conclude that there is a projection $Q'' < P'$ and a partial 
isometry $\nu$ with $\nu \nu^* = Q'$ and $\nu^* \nu = Q''$. Thus $\gamma = [P' -Q'']_{{\bf K_0}(A)}$ in 
${\bf K_0}(A)$. This proves the theorem.
\end{proof}

\begin{cor}
Suppose that $G_1 \supsetneq H \subsetneq G_2$ are countable discrete groups with $H \neq \{ 1 \}$. Suppose 
that $\exists g \in G \overset{def}{=} G_1 \underset{H}{*} G_2$ with $g (H \backslash \{ 1 \}) 
g^{-1} \cap H = \emptyset$. Suppose also that ${\bf K_1}(C^*_r(H)) = 0$. Then 
$${\bf K_0}(C^*_r(G))^+ = \{ \gamma \in {\bf K_0}(C^*_r(G)) | {\bf K_0}(\tau_G)(\gamma) > 0 \} \cup \{ 0 \}.$$
\end{cor}

\begin{proof}
Because of the existence of $\gamma$ we see that condition (\ref{equ:8}) of Theorem \ref{thm:crit} is satisfied.
The existence of $\gamma$ implies also that $H$ is not normal in at least one of the groups $G_1$ or $G_2$. 
Suppose without loss of generality that $H$ is not normal in $G_2$. Then $\Index [G_1:H] \geq 2$ and 
$\Index [G_2:H] \geq 3$ so we can find $g_1 \in G_1 \backslash H$ and $g_2, g_2' \in G_2 \backslash H$ with 
$g_2 g_2'^{-1} \in G_2 \backslash H$. Then condition (\ref{equ:9}) is satisfied with elements $u_1 = \lambda_{g_1}$, $u_2 = \lambda_{g_2}$ and
$u_2' = \lambda_{g_2'}$ and therefore we can apply Theorem \ref{thm:crit}. From the fact that ${\bf
K_1}(C^*_r(H)) = 0$ and Theorem \ref{thm:pimsner} it follows that $\Gamma = {\bf K_0}(C^*_r(G))$. 
This proves the corollary. 
\end{proof}

\begin{remark}
Condition (\ref{equ:9}) is an analogue of the Avitzour condition for the case of reduced amalgamated free products. We
will use it in the next section to prove simplicity and uniqueness of trace.
\end{remark}

\section{Simplicity and Uniqueness of Trace}

In this section we will use Power's idea (\cite{P75}) to obtain a sufficient condition for simplicity and
uniqueness of trace for reduced amalgamated free product $C^*$-algebras. We will make use the following result
(due to Avitzour) and its proof:

\begin{thm}[\cite{A82}] \label{thm:avitzour}
Let $A_1$ and $A_2$ be two unital $C^*$-algebras and $\phi_1$ respectivelly $\phi_2$ states on them with
faithfil GNS-representations. Suppose that there are unitaries $u_i \in A_i$, $i=1,2$ such
that $\phi_1$ and $\phi_2$ are invariant with respect to conjugation by $u_1$ and $u_2$ respectivelly and such
that $\phi_i(u_i) = 0$ for $i=1,2$. Suppose also that there is a unitary $u_2' \in A_2$, such that
$\phi_2(u_2')=0$ and $\phi_2(u_2^* u_2') = 0$. Then: 
\par
(I) $(A,\phi) \overset{def}{=} (A_1, \phi_1) * (A_2, \phi_2)$ is simple. 
\par 
(II) If $\phi$ is invariant with respect to conjugation by $u_2'$ then $\phi$ is the only state on $A$ which
is invariant with respect to conjugation by $u_1, u_2, u_2'$. If $\phi$ is not invariant with respect to 
conjugation by $u_2'$ then there is no state on $A$ which is invariant with respect to conjugation by $u_1,
u_2, u_2'$.
\end{thm}

The proof of Theorem \ref{thm:avitzour} uses a lemma of Choi from \cite{C79}. We will need the following
straightforward generalization of this lemma to the case of Hilbert modules:

\begin{lemma} \label{lemma:choi}
Let $H_1$ and $H_2$ be right Hilbert $B$-modules. Let $u_1, \dots, u_n \in \mathcal{L}(H_1 \oplus H_2)$ be
unitaries such that $u_i^* u_j (H_2) \perp H_2$, whenever $i \neq j$. Suppost that $b \in \mathcal{L}(H_1
\oplus H_2)$ is such that $b(H_1) \perp H_1$. Then $\| \frac{1}{n} \underset{k=1}{\overset{n}{\sum}} u_i^* b u_i \|
\leq 2 \| b \|/\sqrt{n}$.
\end{lemma}

\begin{proof}
First assume that 
$$b = \left[ \begin{array}{cc} 0 & 0 \\ b_1 & b_2 \end{array} \right]  \in \mathcal{L}(H_1
\oplus H_2).$$ 

If $$c = \left[ \begin{array}{cc} c_1 & c_2 \\ 0 & 0 \end{array} \right] \in \mathcal{L}(H_1
\oplus H_2)$$ then for $x \oplus y \in H_1 \oplus H_2$ we have $$ \left[ \begin{array}{cc} c_1 & c_2 \\ b_1 &
b_2 \end{array} \right]  \left( \begin{array}{c} x \\ y \end{array} \right) = \left( \begin{array}{c} c_1x + c_2y \\
b_1x + b_2y \end{array} \right).$$

Then:
$$ \left\| \left( \begin{array}{c} c_1x + c_2y \\ b_1x + b_2y \end{array} \right) \right\|_B^2 = \| \langle (c_1x
+ c_2y) \oplus (b_1x + b_2y), (c_1x + c_2y) \oplus (b_1x + b_2y) \rangle_{H_1 \oplus H_2} \|_B =$$ 
$$=\| \langle c_1x + c_2y, c_1x + c_2y \rangle_{H_1} + \langle b_1x + b_2y, b_1x + b_2y \rangle_{H_2} \|_B 
\leq$$ 
$$\| \langle c_1x + c_2y, c_1x + c_2y \rangle_{H_1} \|_B + \| \langle b_1x + b_2y, b_1x + b_2y 
\rangle_{H_2} \|_B = \| c_1x + c_2y \|_B^2 + \| b_1x + b_2y \|_B^2$$
$$= \left\| \left[ \begin{array}{cc} c_1 & c_2 \\ 0 & 0 \end{array} \right] \left( \begin{array}{c} x \\ y
\end{array} \right) \right\|_B^2 + \left\| \left[ \begin{array}{cc} 0 & 0 \\ b_1 & b_2 \end{array} \right] \left( \begin{array}{c} x \\ y
\end{array} \right) \right\|_B^2.$$

Taking supremum on both sides over all vectors $x \oplus y$ in the unit ball of $H_1 \oplus H_2$ we get
$$ \left\| \left[ \begin{array}{cc} c_1 & c_2 \\ b_1 & b_2 \end{array} \right] \right\|^2 = \| c + b \|^2 \leq
\| c \|^2 + \| b \|^2.$$ 

Now $u_j^* u_i b u_i^* u_j ( H_2 ) \subseteq u_j u_i^* b (H_1) = 0$. So $u_j^* u_i b u_i^* u_j$ has the form
$\left[ \begin{array}{cc} c_1 & c_2 \\ 0 & 0 \end{array} \right]$. Now $\| \underset{i=1}{\overset{n}{\sum}}
u_i b u_i^* \|^2 = \| u_1^* (\underset{i=1}{\overset{n}{\sum}} u_i b u_i^*) u_1 \|^2 = \| b +
\underset{i=2}{\overset{n}{\sum}} u_1^* u_i b u_i^* u_1 \| \leq \| b \|^2 + \|
\underset{i=2}{\overset{n}{\sum}} u_1^* u_i b u_i^* u_1 \|^2 = \| b \|^2 + \|
\underset{i=2}{\overset{n}{\sum}} u_i b u_i^* \|^2.$ It follows by induction that $\|
\underset{i=1}{\overset{n}{\sum}} u_i b u_i^* \|^2 \leq n \| b \|^2$. 
For the general case we represent $$b = \left[ \begin{array}{cc} 0 & b_3 \\ b_1 & b_2 \end{array} \right] = 
\left[ \begin{array}{cc} 0 & 0 \\ b_1 & b_2 \end{array} \right] + 
\left[ \begin{array}{cc} 0 & 0 \\ b_3^* & 0 \end{array} \right]^*.$$ 

Then 
$$\| \underset{i=1}{\overset{n}{\sum}} u_i b u_i^* \| \leq \| \underset{i=1}{\overset{n}{\sum}} u_i \left[
\begin{array}{cc} 0 & 0 \\ b_1 & b_2 \end{array} \right] u_i^* \| + \| \underset{i=1}{\overset{n}{\sum}} u_i
\left[ \begin{array}{cc} 0 & 0 \\ b_3^* & 0 \end{array} \right] u_i^* \| \leq$$ 
$$\sqrt{n} \left\| \left[
\begin{array}{cc} 0 & 0 \\ b_1 & b_2 \end{array} \right] \right\| + \sqrt{n} \left\| \left[ \begin{array}{cc} 0 & 0 \\
b_3^* & 0 \end{array} \right] \right\| \leq 2 \sqrt{n} \| b \|.$$
\end{proof}

Untill the end of the section we will assume that we have unital $C^*$-algebras $A_1$, $A_2$ that contain the 
unital $C^*$-algebra $B$ as a unital $C^*$-subalgebra. We will also assume that we have condiditonal 
expectations $E_i : A_i \rightarrow B$ for $i=1,2$ that have faithful KSGNS-representations (i.e. satisfy 
condition (\ref{equ:1})). We now form the reduced amalgamated free product $(A, E) \overset{def}{=} (A_1, E_1) * (A_2,
E_2)$.
\par
Now we can imitate Avitzour's proof of Theorem \ref{thm:avitzour} and prove the following version for the
amalgamated case:

\begin{prop} \label{prop:avitzour}
Suppose everything is as above and also suppose that there are unitaries $u_1 \in A_1$, $u_2, u_2' \in A_2$
with $E_1(u_1) = 0 = E_2(u_2) = E_2(u_2') = E(u_2 u_2'^*)$. Then if $x \in \Lambda_B^1$ then $0 \in
\overline{\conv} \{ uxu^* | u \in A \text{ is a unitary } \}$. 
\end{prop}

\begin{proof}
We will use the notation from section \ref{S:2} with $I = \{ 1,2 \}$. Let
$W_0 \subset \mathbb{C}(A)$ be the span of all words from $\Lambda_B(A_1^{\circ}, A_2^{\circ})$ that 
either begin with an element $a_1 \in A_1^{\circ}$ or begin with $u_2^* b$ with $b \in B$, or come from $B$. Let
$W_1 \subset \mathbb{C}(A)$ be the span of all words from $\Lambda_B(A_1^{\circ}, A_2^{\circ})$ that 
begin with an element $a_2 \in A_2^{\circ}$ satisfying $E_2(u_2 a_2) = 0$. Denote 
$$H_i \overset{def}{=} \overline{\pi(W_i)\widehat{1_A}} \subset M,\ i=0,1$$
We have $M = H_0 \oplus H_1$ as right Hilbert $B$-module (the orthogonality is with respect to $\langle ., . \rangle_M$). To show this
notice first that $\Span(W_0 \cup W_1)$ is dense in $A$. Therefore $M = H_0 + H_1$. For every word $w_0 \in
W_0$ and every word $w_1 \in W_1$we have $E(w_0^* w_1) = 0$ which is easy to see by considering the three
possible cases for $w_0$. Thus $H_0 \perp H_1$ by linearity.  
\par
We claim that $(u_2^* u_1)^k (H_1) \subseteq H_0$ for $k \neq 0$.
\par
It is enough to prove that $(u_2^* u_1)^k W_1 \subseteq W_0$. 
\par
If $k > 0$ then $(u_2^* u_1)^k W_1$ is spanned by words from $\Lambda_B^{\circ}(A_1^{\circ}, A_2^{\circ})$ 
starting with $u_2^*$. If $k < 0$ then take any word
$w_1 \in W_1$. Then $w_1 = a_2 w_1'$,where $a_2 \in A_2^{\circ}$ satisfies $E(u_2 a_2) = 0$ and $w_1' \in 
\Lambda_B^{\circ}(A_1^{\circ}, A_2^{\circ})$ starts with an element of $A_1^{\circ}$. Then 
$$(u_2^* u_1)^k w_1 = (u_1^* u_2)^{-k} a_2 w_1' = (u_1^* u_2)^{-k-1} u_1^* (u_2 a_2) w_1'$$ 
is a word, starting with $u_1^* \in A_1^{\circ}$. Thus $(u_2^* u_1)^k W_1 \subseteq W_0$. 
\par
Now $u_2'^* x u_2' \in \Lambda_B^2$ and also it is clear that $(u_2'^* x u_2')(W_0) \subseteq W_1$ by considering
the three possibilities for $W_0$ (notice that $E(u_2' u_2^* b) = 0\ \forall b \in B$). Now we can use Lemma \ref{lemma:choi} and get 
$$\| \frac{1}{N} \underset{k=1}{\overset{N}{\sum}} (u_2'^* u_1)^k (u_2'^* x u_2') (u_2'^* u_1)^{-k} \| \leq
\frac{2 \| x \|}{\sqrt{N}}.$$
This implies that $0 \in \overline{\conv} \{ uxu^* | u \in A \text{ is a unitary } \}$.
\end{proof}

We will prove the next technical lemma:

\begin{lemma} \label{lemma:tech}
Suppose that everything is as above and suppose that there are states $\phi_i$ on $A_i$ for $i=1,2$ which are
invariant with respect to $E_i$, $i=1,2$ and satisfy $\phi_1|_B = \phi_2|_B (\overset{def}{=} \phi_B)$, and 
construct $\phi \overset{def}{=} \phi_B \circ E$. 
\par
Suppose that there are two multiplicative sets $1_A \in \tilde{A}_i \subset A_i$ such that
$\Span(\tilde{A}_i)$ is dense in $A_i$, suppose from $a_i \in \tilde{A}_i$ follows $E_i(a_i),\ a_i -
E_i(a_i),\ a_i - \phi_i(a_i) \in \tilde{A}_i$, for 
$i=1,2$, and $B \cap \tilde{A}_1 = B \cap \tilde{A}_2 \overset{def}{=} \tilde{B}$. 
\par
Suppose also that there are two sets of unitaries $\emptyset \neq W_i \subset \tilde{A}_i \cap 
A_i^{\circ}$ such that $(W_i)^* \subset \tilde{A}_i$ for $i=1,2$. Let $u_i \in W_i$, $i=1,2$ and suppose that $\phi$ is invariant with respect to conjugation 
by $u_1$ and $u_2$.
\par
Suppose also that the following condition, similar to condition (\ref{equ:8}), holds: 
\begin{equation} \label{equ:13}
         \begin{cases}
         \forall b_1, \dots, b_l \in \tilde{B}, \text{ with } \phi(b_1) = \dots = \phi(b_l) =0,\ \exists m 
         \in
         \mathbb{N} \text{ and unitaries }  \\
         \nu_{11}, \dots, \nu_{1m}, \nu_{21}, \dots, \nu_{2m} \text{ such that } \nu_{12}, \dots, \nu_{1m} 
         \in
         W_1,\ \nu_{21}, \dots,
         \nu_{2(m-1)} \in W_2, \text{ and: } \\
         \text{either } \nu_{11} \in W_1,\ \nu_{2m} \in W_2 \text{ or } \\ 
         \nu_{11} = 1_{A_1},\ \nu_{2m} \in W_2, \text{ or } \\
         \nu_{11} \in W_1,\ \nu_{2m} = 1_{A_2}, \\ 
         \nu_{11} = 1_{A_1},\ \nu_{2m} = 1_{A_2},\ k \geq 2 \\ 
         \text{ with } E((\nu_{11} \nu_{21} \nu_{12} \cdots \nu_{1m} \nu_{2m}) b_k (\nu_{11} \nu_{21} \nu_{12}
         \cdots \nu_{1m} \nu_{2m})^*) = 0 \text{ for } k=1, \dots, l, \\ 
         (\text{i.e. there are unitaries that conjugate } \tilde{B} \ominus \mathbb{C} 1_B \text{ out of } B) \\ 
         \end{cases}
\end{equation}

Suppose finally that there are unitaries $\omega_1 \in W_1$ and $\omega_2$ with $\omega_2 = 1_A$ or $\omega_2 
\in W_2$, such that $\forall b \in \tilde{B},\ \exists \omega_1^b \in W_1,$ and $\omega_2^b \in W_2$ if $\omega_2 
\in W_2$ or $\omega_2^b = 1$ if $\omega_2 = 1$ with 
$E((\omega_2^b)^* (\omega_1^b)^* b \omega_1 \omega_2) = 0$. 
\par
Then given $x \in \Alg(\tilde{A}_1 \cup \tilde{A}_2)$ with $\phi(x) = 0$ there exist unitaries $\alpha_1, 
\dots, \alpha_s$ with $\alpha_i \in W_{1+(i \mod 2)}$ such that $\alpha_1^* \cdots 
\alpha_s^* x \alpha_s \cdots \alpha_1 \in \Lambda_B^2$. 
\end{lemma}

\begin{proof}
Until the end of this proof we will use the following sets:
\begin{equation*} 
\tilde{\Lambda}_B^1 \overset{def}{=} \ \Span(\underset{k=0}{\overset{\infty}{\bigcup}} (A_1^{\circ} \cap 
\tilde{A}_1) \cdot [(A_2^{\circ} \cap \tilde{A}_2) \cdot (A_1^{\circ} \cap \tilde{A}_1)]^k \subset \mathbb{C}(A),
\end{equation*}
\begin{equation*} 
\tilde{\Lambda}_B^2 \overset{def}{=} \ \Span(\underset{k=0}{\overset{\infty}{\bigcup}} (A_2^{\circ} \cap 
\tilde{A}_2) \cdot [(A_1^{\circ} \cap \tilde{A}_1) \cdot (A_2^{\circ} \cap \tilde{A}_2)]^k ) \subset \mathbb{C}(A),
\end{equation*}
\begin{equation*}
\tilde{\Lambda}_B^{21} \overset{def}{=} \ \Span(\underset{k=1}{\overset{\infty}{\bigcup}} [(A_2^{\circ} \cap 
\tilde{A}_2) \cdot (A_1^{\circ} \cap \tilde{A}_1)]^k ) \subset \mathbb{C}(A),
\end{equation*}
\begin{equation*}
\tilde{\Lambda}_B^{12} \overset{def}{=} \ \Span(\underset{k=1}{\overset{\infty}{\bigcup}} [(A_1^{\circ} \cap 
\tilde{A}_1) \cdot (A_2^{\circ} \cap \tilde{A}_2)]^k ) \subset \mathbb{C}(A).
\end{equation*}

We can write $x = x_B + x_1 + x_2 + x_{12} + x_{21}$, where $x_B \in \Span(\tilde{B})$ with $\phi(x_B) = 0$, 
$x_1 \in \tilde{\Lambda}_B^1$, $x_2 \in \tilde{\Lambda}_B^2$, 
$x_{12} \in \tilde{\Lambda}_B^{12}$ and $x_{21} \in \tilde{\Lambda}_B^{21}$. We will be alternativelly
conjugating $x$ with unitaries from $W_1$ and $W_2$ until we end up with an element of
$\tilde{\Lambda}_B^2$. So at the start we call the words from
$\tilde{\Lambda}_B^1$ "good words". When we conjugate a word 
$w_1 \in \tilde{\Lambda}_B^1$ with $a_2 \in W_2$ we end up with a word 
$a_2 w_1 a_2^* \in \tilde{\Lambda}_B^2$. 
Now we call the words of $\tilde{\Lambda}_B^2$ "good words". If we now take a word $w_2 \in
\tilde{\Lambda}_B^2$ and conjugate it with an element $a_1 \in W_1$ we obtain the word $a_1 w_2 a_1^* \in
\tilde{\Lambda}_B^1$ so we can call the words from $\tilde{\Lambda}_B^1$ "good words". We will show that
proceeding in this way, i.e. alternativelly conjugating $x$ with elements from $W_1$ and $W_2$ we can come to
an element $\alpha_1^* \cdots \alpha_s^* x \alpha_s \cdots \alpha_1 \in \tilde{\Lambda}_B^2$ consisting of a
linear combination of "good words" from $\tilde{\Lambda}_B^2$. This will prove the lemma. 
\par
We have to consider the following $4$ possibilies: 
\par
(i) Take a word $b \in \tilde{B}$. Suppose that the "good words" are in
$\tilde{\Lambda}_B^2$ and we are going to conjugate $b$ with the element $u_1 \in W_1$. Then we obtain 
$$u_1 b u_1^* = E(u_1 b u_1^*) + (u_1 b u_1^* - E(u_1 b u_1^*))$$
for which $(u_1 b u_1^* - E(u_1 b u_1^*)) \in \tilde{A}_1 \cap A_1^{\circ} \subset \tilde{\Lambda}_1$ is a "good word" and the
word $E(u_1 b u_1^*) \in \tilde{B}$ satisfies $\phi(E(u_1 b u_1^*)) = \phi(b)$. Analoguous conclusion can be
drawn if we suppose that the "good words" are in $\tilde{\Lambda}_B^1$ and we are conjugating with the element
$u_2 \in W_2$. 
\par
(ii) Take a word $\gamma_1 \cdots \gamma_{2n} \in \tilde{\Lambda}_B^{12}$ 
($\gamma_i \in A_{1+(i-1 \mod 2)}^{\circ} \cap \tilde{A}_{1+(i-1 \mod 2)}$) 
and conjugate it with a unitary $a_2 \in W_2$ thinking that the "good words" are in $\tilde{\Lambda}_B^1$. 
We get 
\begin{equation*} 
a_2 \gamma_1 \cdots \gamma_{2n-1} \gamma_{2n} a_2^* = a_2 \gamma_1 \cdots \gamma_{2n-1} E(\gamma_{2n} a_2^*)
+ a_2 \gamma_1 \cdots \gamma_{2n-1} (\gamma_{2n} a_2^* - E(\gamma_{2n} a_2^*)).
\end{equation*}

The first word is from $\tilde{\Lambda}_B^{21}$ of the same length $2n$ as the
word $\gamma_1 \cdots \gamma_{2n-1} \gamma_{2n}$ and the second word is from $\tilde{\Lambda}_B^2$, i.e. a
"good word". If we supposed that the good words were in $\tilde{\Lambda}_B^2$
and we were conjugating with a unitary $a_1 \in W_1$ then we would have 
\begin{equation*} 
a_1 \gamma_1 \cdots \gamma_{2n-1} \gamma_{2n} a_2^* = E(a_1 \gamma_1) \gamma_2 \cdots \gamma_{2n-1} \gamma_{2n} a_1^*
+ (a_1 \gamma_1 - E(a_1 \gamma_1)) \gamma_2 \cdots \gamma_{2n-1} \gamma_{2n} a_1^*
\end{equation*}

So again we end up with a word from $\tilde{\Lambda}_B^{21}$ of length $2n$ and a "good word" from 
$\tilde{\Lambda}_B^1$. 
\par
(iii) In a similar way we can treat a word $\gamma_2 \cdots \gamma_{2n+1} \in \tilde{\Lambda}_B^{21}$
($\gamma_i \in A_{1+(i-1 \mod 2)}^{\circ} \cap \tilde{A}_{1+(i-1 \mod 2)}$). If we conjugate with a unitary
$a_2 \in W_2$ knowing that the "good words" are in $\tilde{\Lambda}_B^1$ we end up with 
$$a_2 \gamma_2 \gamma_3 \cdots \gamma_{2n+1} a_2^* = E(a_2 \gamma_2) \gamma_3 \cdots \gamma_{2n+1} a_2^* +
(a_2 \gamma_2 - E(a_2 \gamma_2)) \gamma_3 \cdots \gamma_{2n+1} a_2^*.$$

The first word is from $\tilde{\Lambda}_B^{12}$ and of the same length $2n$ and the second word is from
$\tilde{\Lambda}_B^2$, i.e. a "good word". In the same way if the good words were in $\tilde{\Lambda}_B^2$ and
we were conjugating with a unitary $a_1 \in W_1$ we would obtain 
$$a_1 \gamma_2 \cdots \gamma_{2n} \gamma_{2n+1} a_1^* = a_1 \gamma_2 \cdots \gamma_{2n} E(\gamma_{2n+1} a_1^*)
+ a_1 \gamma_2 \cdots \gamma_{2n} (\gamma_{2n+1} a_1^* - E(\gamma_{2n+1} a_1^*)).$$

The first word is from $\tilde{\Lambda}_B^{12}$ of length $2n$ and the second word is from
$\tilde{\Lambda}_B^1$, i.e. a "good word". 
\par
(iv) Take a word $\gamma_2 \cdots \gamma_{2n} \in \tilde{\Lambda}_B^2$ 
($\gamma_i \in A_{1+(i-1 \mod 2)}^{\circ} \cap \tilde{A}_{1+(i-1 \mod 2)}$). If the "good words" are in 
$\tilde{\Lambda}_B^1$ and if we conjugate this word with the unitary $u_2 \in W_2$, we will get 
$$u_2 \gamma_2 \gamma_3 \cdots \gamma_{2n-1} \gamma_{2n} u_2^* \ \ \ =\ \ \  
E(u_2 \gamma_2) \gamma_3 \cdots \gamma_{2n-1} E( \gamma_{2n} u_2^*) +$$
$$+ (u_2 \gamma_2 - E(u_2 \gamma_2)) \gamma_3 \cdots \gamma_{2n-1} E( \gamma_{2n} u_2^*) + 
E(u_2 \gamma_2) \gamma_3 \cdots \gamma_{2n-1} ( \gamma_{2n} u_2^* - E( \gamma_{2n} u_2^*)) +$$
$$+ (u_2 \gamma_2 - E(u_2 \gamma_2)) \gamma_3 \cdots \gamma_{2n-1} 
( \gamma_{2n} u_2^* - E( \gamma_{2n} u_2^*)).$$

The last word is in $\tilde{\Lambda}_B^2$, so it is a "good word". The second word is in 
$\tilde{\Lambda}_B^{21}$, the third is in $\tilde{\Lambda}_B^{12}$ and the first one is in 
$\tilde{\Lambda}_B^1$ but of length $2n-3$. Since $\phi$ is invariant with respect to conjugation by $u_2$ we
see that $0 = \phi(\gamma_2 \gamma_3 \cdots \gamma_{2n-1} \gamma_{2n}) = \phi(u_2 \gamma_2 \gamma_3 \cdots
\gamma_{2n-1} \gamma_{2n} u_2^*) = \phi(E(u_2 \gamma_2) \gamma_3 \cdots \gamma_{2n-1} E( \gamma_{2n}
u_2^*))$. 
\par
Similarly if we have a word $\gamma_1 \cdots \gamma_{2n-1} \in \tilde{\Lambda}_B^1$ ($\gamma_i \in A_{1+(i-1
\mod 2)}^{\circ} \cap \tilde{A}_{1+(i-1 \mod 2)}$) and if the "good words" are in $\tilde{\Lambda}_B^2$ and if
we conjugate with the unitary $u_1 \in W_1$ we will get 
$$u_1 \gamma_1 \gamma_2 \cdots \gamma_{2n-2} \gamma_{2n-1} u_1^* \ \ \ =\ \ \  
E(u_1 \gamma_1) \gamma_2 \cdots \gamma_{2n-2} E( \gamma_{2n-1} u_1^*) +$$
$$+ (u_1 \gamma_1 - E(u_1 \gamma_1)) \gamma_2 \cdots \gamma_{2n-2} E( \gamma_{2n-1} u_1^*) + 
E(u_1 \gamma_1) \gamma_2 \cdots \gamma_{2n-2} ( \gamma_{2n-1} u_1^* - E( \gamma_{2n-1} u_1^*)) +$$
$$+ (u_1 \gamma_1 - E(u_1 \gamma_1)) \gamma_2 \cdots \gamma_{2n-2} ( \gamma_{2n-1} u_1^* - 
E( \gamma_{2n-1} u_1^*)).$$

Notice that the last word is from $\tilde{\Lambda}_B^1$, so it is a "good word". The second word is from
$\tilde{\Lambda}_B^{12}$ and the third one is from $\tilde{\Lambda}_B^{21}$. The first word is from
$\tilde{\Lambda}_B^2$ but with length $2n-3$. In this case we also can conclude that $0 = \phi(\gamma_1
\gamma_2 \cdots \gamma_{2n-2} \gamma_{2n-1}) = \phi(u_1 \gamma_1 \gamma_2 \cdots
\gamma_{2n-2} \gamma_{2n-1} u_1^*) = \phi(E(u_1 \gamma_1) \gamma_2 \cdots \gamma_{2n-2} E( \gamma_{2n-1}
u_1^*))$. 
\par
From this we can conclude that if we take the word $\gamma_2 \cdots \gamma_{2n} \in \tilde{\Lambda}_B^2$ and
if the "good words" are in $\tilde{\Lambda}_B^1$ then $(u_1 u_2) \gamma_2 \cdots
\gamma_{2n} (u_2^* u_1^*)$ will be the span of some "good words", i.e. belonging to $\tilde{\Lambda}_B^1$,
some words from $\tilde{\Lambda}_B^{21}$, some words from $\tilde{\Lambda}_B^{12}$, and the word from
$\tilde{\Lambda}_B^2$ with length $2n-5$ 
$$E( u_1 E(u_2 \gamma_2) \gamma_3) \gamma_4 \cdots \gamma_{2n-2} E(\gamma_{2n-1} E( \gamma_{2n} u_2^*) u_1^*)
=$$ 
$$= E( u_1 u_2 \gamma_2 \gamma_3) \gamma_4 \cdots \gamma_{2n-2} E(\gamma_{2n-1} \gamma_{2n} u_2^* u_1^*)$$
if $n \geq 3$. Continuing in the same fashion we see that if $l \geq n/2$, $(u_1 u_2)^l \gamma_2 \cdots
\gamma_{2n} (u_2^* u_1^*)$ will be the span of some "good words", i.e. belonging to $\tilde{\Lambda}_B^1$,
some words from $\tilde{\Lambda}_B^{21}$, some words from $\tilde{\Lambda}_B^{12}$, and a word $b \in
\tilde{B}$. Actually it is easy to see that $b = E((u_1 u_2)^l \gamma_2 \cdots \gamma_{2n} 
(u_2^* u_1^*)^l) \in \tilde{B}$ since this is the element which projects onto $B$ under the conditional
expectation. Notice that since $\phi$ is $E$-invariant and also invariant with respect to conjugation by $u_1$ and $u_2$ then $\phi(E((u_1 u_2)^l
\gamma_2 \cdots \gamma_{2n} (u_2^* u_1^*)^l)) = 0$. 
\par
We can now return to the element $x = x_B + x_1 + x_2 + x_{12} + x_{21}$. Set the words from
$\tilde{\Lambda}_B^1$ to be "good words". From the observation above we see that if $l$ is greater that the
length of the longest word appearing in $x_2$, then $(u_1 u_2)^l x_2 (u_2^* u_1^*)^l$ is the span of some
"good words" from $\tilde{\Lambda}_B^1$, some words from $\tilde{\Lambda}_B^{21}$, some words from
$\tilde{\Lambda}_B^{12}$, and some words from $\tilde{B}$, each one of them when evaluated on $\phi$ gives
$0$. But considering cases (i), (ii) and (iii) we can easily conclude that $x' \overset{def}{=} 
(u_1 u_2)^l x (u_2^* u_1^*)^l$ can be written as $x' = x_B' + x_1' + x_{12}' + x_{21}'$ with $x_B'$ being a
span of words from $\tilde{B}$ and satisfying $\phi(x_B') = 0$, $x_1'$ being a span of "good words" from
$\tilde{\Lambda}_B^1$, $x_{12}'$ being a span of words from $\tilde{\Lambda}_B^{12}$ and $x_{21}'$ being a
span of words from $\tilde{\Lambda}_B^{21}$. 
\par
Let $x_B' = \underset{i=1}{\overset{n}{\sum}}\alpha_i b_i$, where $b_i \in \tilde{B}$ and $\alpha_i \in
\mathbb{C}$. $0 = \phi(x_B') = \phi(\underset{i=1}{\overset{n}{\sum}}\alpha_i b_i) =
\underset{i=1}{\overset{n}{\sum}}\alpha_i \phi(b_i)$. Thus $x_B' = \underset{i=1}{\overset{n}{\sum}}\alpha_i 
(b_i - \phi(b_i))$ if we set $b_i' = b_i - \phi(b_i)$ for $i=1, \dots, n$, then $b_i' \in \tilde{B}$ with
$\phi(b_i') = 0 = \phi(u_2 b_i u_2^*)$. So we can apply condition (\ref{equ:13}) to the set of elements  
$\{ b_1', \dots, b_n', E(u_2 b_1' u_2^*), \dots, E(u_2 b_n' u_2^*) \} \subset \tilde{B}$. 
We obtain unitaries $\nu_1, \dots, \nu_m$. Set
$$u = \begin{cases}
      \nu_1 \cdots \nu_m,& \text{ if } \nu_1 \in W_2, \nu_m \in W_2 \\ 
      u_2 \nu_1 \cdots \nu_m,& \text{ if } \nu_1 \in W_1, \nu_m \in W_2, \\ 
      u_2 \nu_1 \cdots \nu_m u_2,& \text{ if } \nu_1 \in W_1, \nu_m \in W_1, \\ 
      \nu_1 \cdots \nu_m u_2,& \text{ if } \nu_1 \in W_2, \nu_m \in W_1.
      \end{cases}$$

Then it is clear that $u^* x_B' u \in \tilde{\Lambda}_B^2$ and the "good words" are in $\tilde{\Lambda}_B^2$.
Then from cases (ii) and (iii) also follows that $x'' \overset{def}{=} u^* x' u$ can be represented as $x'' =
x_2'' + x_{12}'' + x_{21}''$, where $x_2'' \in \tilde{\Lambda}_B^2$ is a span of "good words" and $x_{12}''
\in \tilde{\Lambda}_B^{12}$, $x_{21}'' \in \tilde{\Lambda}_B^{21}$. \\ 
Let $n$ be the number of words from $\tilde{\Lambda}_B^{21}$ and from $\tilde{\Lambda}_B^{12}$ that appear in
the span of $x_{12}'' + x_{21}''$. We will argue by induction on $n$ to conclude the proof of the lemma. \\  
Let $\gamma_1 \cdots \gamma_{2l} \in \tilde{\Lambda}_B^{12}$ ($\gamma_i \in A_{1+(i-1 \mod 2)}^{\circ}
\cap \tilde{A}_{1+(i-1 \mod 2)}$) is a word from the span of $x_{12}''$. (The case $x_{21}''$ is completely
analoguous.) Set 
$$\tilde{u} \overset{def}{=} \begin{cases}
                             \omega_1 \omega_2 (u_1 u_2)^{l-1},& \text{  if } \omega_2 \in W_2, \\ 
                             \omega_1 (u_2 u_1)^{l-1} u_2,& \text{  if } \omega_2 = 1_A.
                             \end{cases}$$

Let's observe first that if $\alpha_1 \cdots \alpha_{2l}, \beta_1 \cdots \beta_{2l} \in 
\tilde{\Lambda}_B^{12}$, then we can write 
$$E(\beta_{2l}^* \cdots \beta_2^* \beta_1^* \alpha_1 \alpha_2 \cdots \alpha_{2l}) = E(\beta_{2l}^* \cdots
\beta_2^* E(\beta_1^* \alpha_1) \alpha_2 \cdots \alpha_{2l}) + $$ 
$$+ E(\beta_{2l}^* \cdots \beta_2^* (\beta_1^*
\alpha_1 - E(\beta_1^* \alpha_1)) \alpha_2 \cdots \alpha_{2l}) = E(\beta_{2l}^* \cdots \beta_2^* 
E(\beta_1^* \alpha_1) \alpha_2 \cdots \alpha_{2l}).$$ 

It follows by induction that $E(\beta_{2l}^* \cdots
\beta_2^* \beta_1^* \alpha_1 \alpha_2 \cdots \alpha_{2l}) \in \tilde{B}$. Also from 
$$\beta_{2l}^* \cdots \beta_2^* \beta_1^* \alpha_1 \alpha_2 \cdots \alpha_{2l} = \beta_{2l}^* \cdots
\beta_2^* E(\beta_1^* \alpha_1) \alpha_2 \cdots \alpha_{2l} + $$ 
$$+ \beta_{2l}^* \cdots \beta_2^* (\beta_1^*
\alpha_1 - E(\beta_1^* \alpha_1)) \alpha_2 \cdots \alpha_{2l} = \beta_{2l}^* \cdots \beta_2^* 
E(\beta_1^* \alpha_1) \alpha_2 \cdots \alpha_{2l}$$                          
again by induction follows that $\beta_{2l}^* \cdots \beta_2^* \beta_1^* \alpha_1 \alpha_2 \cdots \alpha_{2l}$
is the span of words from $\tilde{\Lambda}_B^2$ plus the word $E(\beta_{2l}^* \cdots \beta_2^* \beta_1^* 
\alpha_1 \alpha_2 \cdots \alpha_{2l}) \in \tilde{B}$. 
\par
All this implies that $\tilde{u}^* \gamma_1 \cdots \gamma_{2l} \tilde{u}$ is a span 
of "good words" from $\tilde{\Lambda}_B^2$ and the word $E(\tilde{u}^* \gamma_1 \cdots \gamma_{2l}) \tilde{u}
\in \tilde{\Lambda}_B^{12}$. Set $\tilde{b} \overset{def}{=} 
E(\tilde{u}^* \gamma_1 \cdots \gamma_{2l}) \in \tilde{B}$ (see the observation above). 
Now we choose unitaries $\omega_1^{\tilde{b}}, \omega_2^{\tilde{b}}$ as in the 
statement of the lemma. We have 
$$(\omega_2^{\tilde{b}})^* (\omega_1^{\tilde{b}})^* E(\tilde{u}^* \gamma_1
\cdots \gamma_{2l}) \tilde{u} \omega_1^{\tilde{b}} \omega_2^{\tilde{b}} = $$
$$= \begin{cases} 
    (\omega_2^{\tilde{b}})^* (\omega_1^{\tilde{b}})^* E(\tilde{u}^* \gamma_1
\cdots \gamma_{2l}) \omega_1 \omega_2 (u_1 u_2)^{l-1} \omega_1^{\tilde{b}} \omega_2^{\tilde{b}},& \text{  if } 
\omega_2 \in W_2, \\ 
    (\omega_1^{\tilde{b}})^* E(\tilde{u}^* \gamma_1
\cdots \gamma_{2l}) \omega_1 (u_2 u_1)^{l-1} u_2 \omega_1^{\tilde{b}},& \text{  if } \omega_2 = 1_A.
     \end{cases}$$
     
From this and from the choice of 
$\omega_1^{\tilde{b}}, \omega_2^{\tilde{b}}$ (and from case (i)) it is clear
that 
$$(\omega_2^{\tilde{b}})^* (\omega_1^{\tilde{b}})^* E(\tilde{u}^* \gamma_1
\cdots \gamma_{2l}) \tilde{u} \omega_1^{\tilde{b}} \omega_2^{\tilde{b}}$$ 
is a span of
"good words". 
\par
Since by cases (ii) and (iii) follows that when we alternatively conjugate words from
$\tilde{\Lambda}_B^{12}$ and from $\tilde{\Lambda}_B^{21}$ by unitaries from $W_1$ and $W_2$ the number of
such words doesn't increase, and since we managed to conjugate the word $\gamma_1 \cdots \gamma_{2l}$ to a
span of "good words", the induction on $n$ is concluded. 
\par
This proves the lemma. 
\end{proof}

Combining Proposition \ref{prop:avitzour} and Lemma \ref{lemma:tech} we obtain the following 

\begin{thm} \label{thm:dualcrit}
Assume that we have unital $C^*$-algebras $A_1$, $A_2$ that contain the unital $C^*$-algebra $B$ as a unital 
$C^*$-subalgebra. Also assume that there are condidional expectations $E_i : A_i \rightarrow B$ for $i=1,2$ 
that have faithful KSGNS-rapresentations (i.e. satisfy condition (\ref{equ:1})) and form the reduced amalgamated free 
product $(A, E) \overset{def}{=} (A_1, E_1) * (A_2, E_2)$. 
\par
Suppose that there are states $\phi_i$ on $A_i$ for $i=1,2$ which are invariant with respect to $E_i$,
$i=1,2$ and satisfy $\phi_1|_B = \phi_2|_B (\overset{def}{=} \phi_B)$. Construct 
$\phi \overset{def}{=} \phi_B \circ E$. \\ 
Assume that there are unitaries $u_1 \in A_1$, $u_2, u_2' \in A_2$ with $E_1(u_1) = 0 = E_2(u_2) = E_2(u_2') 
= E(u_2 u_2'^*)$. (Or assume that there are unitaries $u_1, u_1' \in A_1^{\circ}$, $u_2 \in A_2^{\circ}$ with
$E(u_1 u_1'^*) = 0$.) 
\par
Suppose that there are two multiplicative sets $1_A \in \tilde{A}_i \subset A_i$ such that
$\Span(\tilde{A}_i)$ is dense in $A_i$, suppose from $a_i \in \tilde{A}_i$ follows $E_i(a_i),\ a_i -
E_i(a_i),\ a_i - \phi_i(a_i) \in \tilde{A}_i$, for 
$i=1,2$, and $B \cap \tilde{A}_1 = B \cap \tilde{A}_2 \overset{def}{=} \tilde{B}$. 
\par
Suppose also that there are two sets of unitaries $\emptyset \neq W_i \subset \tilde{A}_i \cap 
A_i^{\circ}$ such that $(W_i)^* \subset \tilde{A}_i$ for $i=1,2$. Let $v_i \in W_i$, $i=1,2$ and suppose that $\phi$ is invariant with respect to conjugation 
by $v_1$ and $v_2$. 
\par
Suppose that condition (\ref{equ:13}) holds, namely:
\begin{equation*} 
         \begin{cases}
         \forall b_1, \dots, b_l \in \tilde{B}, \text{ with } \phi(b_1) = \dots = \phi(b_l) =0,\ \exists m 
         \in
         \mathbb{N} \text{ and unitaries }  \\
         \nu_{11}, \dots, \nu_{1m}, \nu_{21}, \dots, \nu_{2m} \text{ such that } \nu_{12}, \dots, \nu_{1m} 
         \in
         W_1,\ \nu_{21}, \dots,
         \nu_{2(m-1)} \in W_2, \text{ and: } \\
         \text{either } \nu_{11} \in W_1,\ \nu_{2m} \in W_2 \text{ or } \\ 
         \nu_{11} = 1_{A_1},\ \nu_{2m} \in W_2, \text{ or } \\
         \nu_{11} \in W_1,\ \nu_{2m} = 1_{A_2}, \\ 
         \nu_{11} = 1_{A_1},\ \nu_{2m} = 1_{A_2},\ k \geq 2 \\ 
         \text{ with } E((\nu_{11} \nu_{21} \nu_{12} \cdots \nu_{1m} \nu_{2m}) b_k (\nu_{11} \nu_{21} \nu_{12}
         \cdots \nu_{1m} \nu_{2m})^*) = 0 \text{ for } k=1, \dots, l, \\ 
         (\text{i.e. there are unitaries that conjugate } \tilde{B} \ominus \mathbb{C} 1_B \text{ out of } B) \\ 
         \end{cases}
\end{equation*} 

Suppose finally that there are unitaries $\omega_1 \in W_1$ and $\omega_2$ with $\omega_2 = 1_A$ or $\omega_2 
\in W_2$, such that $\forall b \in \tilde{B},\ \exists \omega_1^b \in W_1,$ and $\omega_2^b \in W_2$ if $\omega_2 
\in W_2$ or $\omega_2^b = 1$ if $\omega_2 = 1$ with 
$E((\omega_2^b)^* (\omega_1^b)^* b \omega_1 \omega_2) = 0$. 
\par
Then: 
\par
(1) If $\phi_B$ has a faithful GNS-representation then $A$ is simple. 
\par
(2) If $\phi$ is invariant with respect to conjugation by $u_1, u_2, u_2'$ (or by $u_1, u_1', u_2$) and all the unitaries from $W_1$
and $W_2$, then $\phi$ is the only tracial state on $A$, invariant with respect to conjugation by all those
unitaries.
\end{thm}

\begin{proof}
(1) Suppose $I  \neq 0$ is an ideal of $A$. Notice that $\Alg(\tilde{A}_1 \cup \tilde{A}_2)$ is dense in $A$. 
Take a nonzero element $x \in I$. 
Because $E$ has
a faithful KSGNS-representation it satisfies condition (\ref{equ:1}), i.e. $\exists y \in A$ such that $b
\overset{def}{=} E(y^* x^* x y) \neq 0$. Notice that $b^* = b$. Since $\phi_B$ has a faithful GNS-representation 
we can find $b' \in B$ such that $\phi_B((b')^*  b b') \neq 0$. Then 
$$\phi((b')^* y^* x^* x y b') = \phi(E((b')^* y^* x^* x y b')) = \phi((b')^* E(y^* x^* x y) b') =
\phi((b')^* b b') \neq 0.$$

Then $c \overset{def}{=} \phi((b')^* b b')^{-1} (b')^* y^* x^* x y b' $ ($\in I$) is self-adjoined and satisfies
$\phi(c) = 1$. Find $a \in \Alg(\tilde{A}_1 \cup \tilde{A}_2)$ such that $\| a - c \| < \epsilon$. From 
Lemma \ref{lemma:tech} it follows that we can find unitaries $\alpha_1, \dots, \alpha_m \in W_1 \cup W_2$ 
such that $(\alpha_1 \cdots \alpha_m)^* (a - \phi(a)1_A) (\alpha_1 \cdots \alpha_m) \in \Lambda_B^1$.
Then it follows from Proposition \ref{prop:avitzour} that we can find unitaries $U_1, \dots U_N \in A$ that
are constructed from $u_1, u_2, u_2'$ and the unitaries from $W_1 \cup W_2$ and are such that 
$$\| \underset{i=1}{\overset{N}{\sum}} \frac{1}{N} U_i^* (\alpha_1 \cdots \alpha_m)^* (a - \phi(a)1_A) 
(\alpha_1 \cdots \alpha_m) U_i \| < \epsilon.$$ 

Then 
$$\| \underset{i=1}{\overset{N}{\sum}} \frac{1}{N} U_i^* (\alpha_1 \cdots \alpha_m)^* (a - \phi(a)1_A - c +
1_A) (\alpha_1 \cdots \alpha_m) U_i \| \leq$$ 
$$\underset{i=1}{\overset{N}{\sum}} \frac{1}{N} \| U_i^* (\alpha_1 \cdots \alpha_m)^* 
(a - \phi(a)1_A - c + 1_A) (\alpha_1 \cdots \alpha_m) U_i \| = $$
$$ = \underset{i=1}{\overset{N}{\sum}} \frac{1}{N} \|  a - \phi(a)1_A - c + 1_A \| = \| a - \phi(a)1_A - c + 1_A
\| =$$
$$ = \| (a-c) - \phi(a-c) \| \leq \|a-c\| + \|a-c\| < 2 \epsilon.$$

Therefore $\| \underset{i=1}{\overset{N}{\sum}} \frac{1}{N} U_i^* (\alpha_1 \cdots \alpha_m)^* (c - 1_A) 
(\alpha_1 \cdots \alpha_m) U_i \| < 3 \epsilon$. Set \\ 
$d \overset{def}{=}  \underset{i=1}{\overset{N}{\sum}} 
\frac{1}{N} U_i^* (\alpha_1 \cdots \alpha_m)^* c 
(\alpha_1 \cdots \alpha_m) U_i$ ($\in I$).
\par
Thus $\| d - 1_A \| < 3 \epsilon$. Then if we take $\epsilon < \frac{1}{3}$ it would follow that $d$ is 
invertible, and therefore $I = A$. 
\par
(2) Take $0 \neq x \in A$. Then if we argue as in case (1) we can find unitaries $U_1, \dots, U_N \in 
\overline{conv} \{ u | u \text{ is a product of unitaries from } W_1 \cup W_2 \cup \{ u_1, u_2, u_2' \} \}$ 
with $$\| \underset{i=1}{\overset{N}{\sum}} \frac{1}{N} U_i^* (x - \phi(x)1_A) U_i \| < 3 \epsilon.$$ 
If we take a state $\phi'$ such that $\phi$ and $\phi'$ are invariant with respect to conjugation by $u_1,
u_2, u_2'$ and by all unitaries from $W_1 \cup W_2$ then we will have 
$$3 \epsilon > | \phi'(\underset{i=1}{\overset{N}{\sum}} \frac{1}{N} U_i^* (x - \phi(x)1_A) U_i) | = |
\underset{i=1}{\overset{N}{\sum}} \frac{1}{N} \phi'(U_i^* x U_i) - \phi(x) | = |
\underset{i=1}{\overset{N}{\sum}} \frac{1}{N} \phi'(x) - \phi(x) | =$$ 
$$|\phi'(x) - \phi(x) |.$$

Since this is true for any $\epsilon > 0$ it follows that $\phi' \equiv \phi$.
\end{proof}

Although the statement of Theorem \ref{thm:dualcrit} looks complicated some 
applications can be given. The  next proposition is a slight generalization of the de la 
Harpe's result from \cite{dlH85}.

\begin{cor}
Suppose that $G_1 \supsetneq H \subsetneq G_2$ are discrete groups and suppose that $H \neq \{ 1 \}$. Denote $G 
\overset{def}{=} G_1 \underset{H}{*} G_2$. Suppose that for any finitely many $h_1, \dots, h_m \in H 
\backslash \{ 1 \}$ there is $g \in G$ with $g^{-1} h_i g \notin H$ for all $i=1, \dots, m$.
Then $C^*_r(G)$ is simple with a unique trace.
\end{cor}

\begin{proof}
Set $A_i = C^*_r(G_i)$, $i=1,2$, $B = C^*_r(H)$ and $A = C^*_r(G)$. Clearly $H$ is not normal in at least one
of the groups $G_1$ or $G_2$. Without loss of generality suppose that $H$ is not normal in
$G_1$. Then there are $g_1, g_1' \in G_1 \backslash H$ and $g_2 \in G_2 \backslash H$ with $g_1 (g_1')^{-1}
\in G_1 \backslash H$. Then set $u_1 = \lambda(g_1)$, $u_1' = \lambda(g_1')$, $u_2 = \lambda(g_2)$. We take
$\tilde{A}_i = \{ \lambda(c_i) | c_i \in G_i \}$, $i = 1,2$, $\tilde{B} = \{ \lambda(h) | h \in H \}$. 
Also $W_i = \tilde{A}_i \backslash \tilde{B}$ for $i=1,2$. Condition (\ref{equ:13}) is satisfied since for finitely many
elements from $H \backslash \{ 1 \}$ we can find an element from $G$ that conjugates them away from $H$.
Finally for the last condition of Theorem \ref{thm:dualcrit} we can set $\omega_1 = u_1$, $\omega_2 = 1$ and
for $\lambda(h) \in \tilde{B}$ we set $\omega_1^{\lambda(h)} = h u_1'$. Thus all
requirements of Theorem \ref{thm:dualcrit} are met and this finishes the proof.
\end{proof}

We give also an application to HNN extensions of discrete groups. We will use the notion of reduced
HNN extensions for $C^*$-algebras introduced by Ueda in \cite{U06}. We will use the following settings: 
\par
Let $\{ 1 \} \subsetneq H \subset G$ be countable discrete groups and let $\tilde{\theta} : H \to G$ be an 
injective group homomorphism. Thus we have that $C^*_r(H) \subset C^*_r(G)$ and that we have a well defined
injective $*$-homomorphism $\theta : C^*_r(H) \to C^*_r(G)$. By $E^G_H : C^*_r(G) \to C^*_r(H)$ and 
$E^G_{\tilde{\theta}(H)} : C^*_r(G) \to C^*_r(\theta(C^*_r(H))$ we will denote the canonical conditional expectations.
By $\tau_G$ we will denote the canonical trace on $C^*_r(G)$. Let $A_1 = C^*_r(G) \otimes M_2(\mathbb{C})$,
$A_2 = C^*_r(H) \otimes M_2(\mathbb{C})$ and $B = C^*_r(H) \oplus C^*_r(H)$. Define the inclusion maps $i_1 :
B \to A_1$ and $i_2 : B \to A_2$ as
$$i_1(b_1 \oplus b_2) = \left[  \begin{array}{cc} b_1 & 0 \\ 0 & \theta(b_2) \end{array} \right], \qquad 
i_2(b_1 \oplus b_2) = \left[  \begin{array}{cc} b_1 & 0 \\ 0 & b_2 \end{array} \right]$$
and define the conditional expectations $E_1 : A_1 \to B$ and $E_2 :A_2 \to B$ as 
$$E_1 = \left[  \begin{array}{cc} E^G_H & 0 \\ 0 & E^G_{\tilde{\theta}(H)} \end{array} \right], \qquad 
E_2 = \left[  \begin{array}{cc} \id & 0 \\ 0 & \id \end{array} \right].$$

Then let $$(A, E) = (A_1, E_1) * (A_2, E_2)$$ be the reduced amalgamated free product of $(A_1, E_1)$ and
$(A_2, E_2)$ and let $$(\mathcal{A}, E^{\mathcal{A}}_{C^*_r(G)}, u(\theta)) = (C^*_r(G), E^G_H) \bigstar_{C^*_r(H)}
(\theta, E^G_{\tilde{\theta}(H)})$$ be the reduced HNN extension of $C^*_r(G)$ by $\theta$ as in
\cite{U06}. Also let $i_B : B \to A$ be the canonical inclusion. 
\par
From \cite[Proposition 2.2]{U06} follows that $A$ is isomorphic to $\mathcal{A} \otimes M_2(\mathbb{C})$. Therefore the
questions of simplicity and uniqueness of trace for $A$ and for $\mathcal{A}$ are equivalent. The following 
corollary of Theorem  \ref{thm:dualcrit} is true:

\begin{cor}
In the above settings suppose that $H \subsetneq G$ and $\tilde{\theta}(H) \subsetneq G$. Suppose also that
$\forall h \in H \backslash \{ 1 \}$, $\exists n_h \in \mathbb{N}$, such that 
$\tilde{\theta}^{n_h-1}(h) \in H$ and $\tilde{\theta}^{n_h}(h) \notin H$. 
Then $A$ (and therefore $\mathcal{A}$ also) is simple with a unique trace.
\end{cor}

\begin{proof}
We will show that all the conditions of Theorem \ref{thm:dualcrit} are met.
\par
First the canonical traces $\tau_i$ on $A_i$, $i=1,2$ satisfy $\tau_i \circ E_i = \tau_i$ for $i=1,2$ and
$\tau_1|_B = \tau_2|_B$ ($\overset{def}{=} \tau_B$). We have $\tau = \tau_B \circ E$. 
\par
Define $$\tilde{A}_1 = \Span(\{ \lambda(g) \otimes e_{ij} | g \in G,\ 1 \leq i,j \leq 2 \})$$ and 
$$\tilde{A}_2 = \Span(\{ \lambda(h) \otimes e_{ij} | h \in H,\ 1 \leq i,j \leq 2 \}),$$ where $e_{ij}$ for $1 \leq i,j \leq 2$ are the
matrix units for $M_2(\mathbb{C})$. Then we have $\tilde{A}_1 \cap B = \tilde{A}_2 \cap B$($\overset{def}{=}
\tilde{B}$). It is also clear that $a_i \in \tilde{A}_i$ implies $E(a_i),\ a_i - E(a_i), a_i - \tau_i(a_i) \in
\tilde{A}_i$ for $i=1,2$. 
\par
Choose $\bar{g}_1 \in G \backslash H,\ \bar{g}_2 \in G \backslash \tilde{\theta}(H)$. 
\par
Define the following unitaries from $A_1 \cap \tilde{A}_1$: 
$$u_1 = \left[ \begin{array}{cc} 0 & 1 \\ 1 & 0 \end{array} \right], \qquad 
u_1' = \left[  \begin{array}{cc} \lambda(\bar{g}_1) & 0 \\ 0 & \lambda(\bar{g}_2) \end{array} \right], \qquad 
u_1'' =  \frac{1}{\sqrt{2}}\left[ \begin{array}{cc} -\lambda(\bar{g}_1) & \lambda(\bar{g}_1) \\
\lambda(\bar{g}_2) & \lambda(\bar{g}_2) \end{array} \right],$$ and the following unitary from $A_2 \cap
\tilde{A}_2$:
$$u_2 = \left[ \begin{array}{cc} 0 & 1 \\ 1 & 0 \end{array} \right].$$

Set $W_1 = \{ u_1, u_1', u_1'' \}$, $W_2 = \{ u_2 \}$.  
\par
Set $\omega_1 = u_1$, $\omega_2 = 1_{A_2}$ and for every $b = b_1 \oplus b_2 \in \tilde{B}$ set 
$\omega_1^b = u_1'$. Then 
$$E((u_1')^* (b_1 \oplus b_2) u_1) = 
E( \left[  \begin{array}{cc} \lambda(\bar{g}_1^{-1}) & 0 \\ 0 & \lambda(\bar{g}_2^{-1}) \end{array} \right] 
\left[  \begin{array}{cc} b_1 & 0 \\ 0 & \theta(b_2) \end{array} \right] 
\left[ \begin{array}{cc} 0 & 1 \\ 1 & 0 \end{array} \right]) = $$
$$= E(\left[ \begin{array}{cc} 0 & \lambda(\bar{g}_1^{-1})b_1 \\ \lambda(\bar{g}_1^{-1})\theta(b_2) & 0 \end{array} \right])
=0.$$

It remains to check that condition (\ref{equ:13}) holds.  
\par
For an element $b = b_1 \oplus b_2 \in B$ it is easy to see that 
$$u_2^* u_1^* b u_1 u_2 = E(u_2^* E(u_1^* b u_1) u_2) + u_2^*(u_1^*bu_1 - E(u_1^*bu_1))u_2$$ and
that 
$$i_B^{-1} \circ E(u_2^* u_1^* b u_1 u_2) = 
      \begin{cases} 
      \theta^{-1}(b_1) \oplus \theta(b_2),&\ \ \text{ if } b_1 \in \theta(C^*_r(H)),\ b_2 \in C^*_r(H), \\ 
      \theta^{-1}(b_1) \oplus 0,&\ \ \text{ if } b_1 \in \theta(C^*_r(H)),\ b_2 \notin C^*_r(H), \\ 
      0 \oplus \theta(b_2),&\ \ \text{ if } b_1 \notin \theta(C^*_r(H)),\ b_2 \in C^*_r(H), \\ 
      0 \oplus 0,&\ \ \text{ if } b_1 \notin \theta(C^*_r(H)),\ b_2 \notin C^*_r(H).
      \end{cases}$$

Using induction one can show that for any $n \in \mathbb{N}$ we have 
$$(u_2^* u_1^*)^n b (u_1 u_2)^{-n} - E((u_2^* u_1^*)^{-n} b (u_1 u_2)^n) \in \Lambda_B^2.$$ 

Let $\hat{\theta}$ be the linear map which extends $\theta$ to $C^*_r(G)$ by $\hat{\theta}(\lambda(g)) = 0$ 
for $g \in G \backslash H$. Also let $\hat{\theta^{-1}}$ be the linear map which extends $\theta^{-1}$ to 
$C^*_r(G)$ by $\hat{\theta^{-1}}(\lambda(g)) = 0$ for $g \in G \backslash \tilde{\theta}(H)$. Then:
$$i_B^{-1} \circ E((u_2^* u_1^*)^{-n} b (u_1 u_2)^n) =$$
$$=
      \begin{cases} 
      \theta^{-n}(b_1) \oplus \theta^n(b_2),& \text{ if } b_1 \in (\hat{\theta})^n(C^*_r(H)),\ 
      b_2 \in (\hat{\theta^{-1}})^{n-1}(C^*_r(H)), \\ 
      \theta^{-n}(b_1) \oplus 0,& \text{ if } b_1 \in (\hat{\theta})^n(C^*_r(H)),\ b_2 \notin 
      (\hat{\theta^{-1}})^{n-1}(C^*_r(H)), \\ 
      0 \oplus \theta^n(b_2),& \text{ if } b_1 \notin (\hat{\theta})^n(C^*_r(H)),\ b_2 \in 
      (\hat{\theta^{-1}})^{n-1}(C^*_r(H)), \\ 
      0 \oplus 0,& \text{ if } b_1 \notin (\hat{\theta})^n(C^*_r(H)),\ b_2 \notin 
      (\hat{\theta^{-1}})^{n-1}(C^*_r(H)).
      \end{cases}$$

If we set $c_1 = \lambda(\bar{g}_1^{-1}) (\hat{\theta^{-1}})^{n}(b_1) \lambda(\bar{g}_1)$ and 
$c_2 = \lambda(\bar{g}_2^{-1}) (\hat{\theta})^{n+1}(b_2) \lambda(\bar{g}_2)$ the we will have
$$i_B^{-1} \circ E(u_2^* (u_1')^* (u_2^* u_1^*)^n b (u_1 u_2)^n u_1' u_2) =$$
$$=
      \begin{cases} 
      \theta^{-1}(c_2) \oplus c_1,&\ \text{ if } c_2 \in \theta(C^*_r(H)),\ c_1 \in C^*_r(H), \\ 
      \theta^{-1}(c_2) \oplus 0,&\ \text{ if } c_2 \in \theta(C^*_r(H)),\ c_1 \notin C^*_r(H), \\ 
      0 \oplus c_1,&\ \text{ if } c_2 \notin \theta(C^*_r(H)),\ c_1 \in C^*_r(H), \\ 
      0 \oplus 0,&\ \text{ if } c_2 \notin \theta(C^*_r(H)),\ c_1 \notin C^*_r(H).
      \end{cases}$$

Now take elements $\tilde{b}_1, \dots, \tilde{b}_l \in \tilde{B}$ with $\tau_B(\tilde{b}_1) = \dots =
\tau_B(\tilde{b}_l) = 0$. We can write $\tilde{b}_k = \alpha_k + b_{k1} \oplus -\alpha_k + b_{k2}$ for each
$k=1, \dots, l$ with $b_{kj} \in \Span(\{ \lambda(h) | h \in H \backslash \{ 1 \} \})$. \\ 
Clearly from the statement of the corollary follows that there exists an $N \in \mathbb{N}$ with 
$E^G_H(\hat{\theta}^N(b_{k2})) = 0$ for each $k=1, \dots, l$.
Therefore for each $k = 1, \dots, l$ we have 
$$i_B^{-1} \circ E(u_2^* (u_1')^* (u_2^* u_1^*)^{-N} \tilde{b}_k (u_1 u_2)^N u_1' u_2) =$$
$$=
     \begin{cases}
     \alpha_k \oplus -\alpha_k + c_k,&\ \text{ if } c_k \in C^*_r(H), \\ 
     \alpha_k \oplus -\alpha_k,&\ \text{ if } c_k \notin C^*_r(H),
     \end{cases}$$ 
where $c_k = \lambda(\bar{g}_1^{-1}) (\hat{\theta^{-1}})^{N}(b_{k1}) \lambda(\bar{g}_1)$, $k=1, \dots l$. Now we can find
an $M \in \mathbb{N}$ such that $(\hat{\theta})^M(c_k) = 0$ for all $k=1, \cdots, l$. Then for all $k=1, \dots, l$ we 
have
$$i_B^{-1} \circ E((u_2^* u_1^*)^{-M}  u_2^* (u_1')^* (u_2^* u_1^*)^{-N} \tilde{b}_k (u_1 u_2)^N u_1' u_2 (u_1 u_2)^M) =
\alpha_k \oplus -\alpha_k.$$

Finally for all $k=1, \dots, l$ 
$$i_B^{-1} \circ E((u_1'')^*(u_2^* u_1^*)^{-M}  u_2^* (u_1')^* (u_2^* u_1^*)^{-N} \tilde{b}_k (u_1 u_2)^N u_1' u_2 (u_1
u_2)^M u_1'') =0.$$

This proves that condition (\ref{equ:13}) holds and thus we can apply Theorem \ref{thm:dualcrit}. 
\par
This proves the Corollary.
\end{proof}

\begin{remark}
By symmetry it is clear that in the corollary the assumption \\
$''$ $\forall h \in H \backslash \{ 1 \}$, $\exists n_h \in \mathbb{N}$, such that 
$\tilde{\theta}^{n_h-1}(h) \in H$ and 
$\tilde{\theta}^{n_h}(h) \notin H$ $''$ \\
can be replaced by the assumption \\ 
$''$ $\forall h \in H \backslash \{ 1 \}$, $\exists n_h \in \mathbb{N}$, such that 
$\tilde{\theta}^{-n_h+1}(h) \in \tilde{\theta}(H)$ and 
$\tilde{\theta}^{-n_h}(h) \notin \tilde{\theta}(H)$ $''$.
\end{remark}

Examples of HNN extensions of discrete groups which satisfy the assumption of this corollary (and whose
reduced $C^*$-algebras are simple with a unique trace) are the Baumslag-Solitar groups $BS(n,m) \overset{def}{=} 
\langle a,b\ |\ b^{-1} a^m b = a^n \rangle$ for $|n| \neq |m|$ and $|n|, |m| \geq 2$. 
\par
Somewhat related result is the ICC property. It was proved by Stalder in \cite{S06}
that $BS(m,n)$ is an ICC group if and only if $|n| \neq |m|$. \\
\par 
Prof. Ueda pointed out to me that our result on the $C^*$-simplicity of $BS(m,n)$ is sharp: 
\par
In the case $m=1$ (or
$n=1$) it is known that those groups are solvable. $B(1,1)$ is abelian. For $|n| > 1$ $BS(1,n) =\langle a,b\ |\
b^{-1} a b = a^n \rangle$. It is not difficult to see that all the elements of $BS(1,n)$ can be written in the
form $b^i a^k b^{-j}$, where $i,j \geq 0$ and if $ij > 0$ then $n \nmid k$. Then $b^i a^k b^{-j} \mapsto i-j$ is
a well defined group homomorphism $h : BS(1,n) \to \mathbb{Z}$. Then one can check that $\ker(h) = \langle 
b^i a^k b^{-i},\ i \geq 0,\ k \in \mathbb{Z}\ |\ b^{i+1} a^{nk} b^{-i-1} = b^i a^k b^{-i} \rangle$ and it is isomorphic to the
additive group of the $n$-adic numbers. This shows that $BS(1,n)$ is meta-abelian (extension of an abelian by an 
abelian) group and trerefore solvable. It is also know that extension of an amenable group by an amenable group
is amenable group and therefore the solvable groups are amenable. 
\par
If $G$ is a locally compact discrete group then we have (by definition) the one-dimensional representation of the 
full $C^*$-algebra of $G$ $\pi : C^*(G) \to \mathbb{C}$ given on the generators of $G$ by $\pi(f_g) = 1$ for all 
$g \in G$ ($f_g: G \to \mathbb{C}$ is given by $f_g(h) = \delta_{gh},\ h \in G$). If $|G| > 1$ $\ker(\pi)$ is a
nontrivial ideal in $C^*(G)$. Obviously $\pi$ is a tracial state. If $|G| > 1$ then $1 = \pi(f_g) \neq \tau_G(f_g)
= 0$ for $\forall 1 \neq g \in G$, where $\tau_G$ is the canonical trace on $C^*(G)$. Therefore if $|G| > 1$
then $C^*(G)$ has more that one trace. All this shows that if $G$ is an amenable locally compact discrete group 
and if $|G| > 1$ then $C^*_r(G)$ ($= C^*(G)$) is not simple and has more than one trace. Therefore
$C^*_r(BS(1,n))$ are not simple and each one has more than one trace for each $n \in \mathbb{Z}$. 
\par
Finally if $m = n$ and $|n|, |m| \geq 2$ then $BS(n,n)$ has a nontrivial center ($a^n$ is in the center of
$BS(n,n)$). If $m = - n$ then $C^*_r(BS(-n,n))$ has a nontrivial center ($\lambda(a^n) + \lambda(a^{-n})$ is in 
the center of $C^*_r(BS(-n,n))$). In both cases $C^*_r(BS(m,n)) \cong A \otimes C(X)$ for some $C^*$-algebra $A$
and some compact Hausdorff space $X$ ($|X| > 1$). If $x \in X$ then $I = \langle a \otimes f\ |\ a \in A,\ f \in
C(X),\ f(x) = 0 \rangle_{A \otimes C(X)}$ is a nontrivial ideal of $A \otimes C(X)$. Also if we call $\tau_A$ the
restriction of the canonical trace on $C^*_r(BS(-n,n))$ to $A \otimes 1_X$ then if $x, y \in X$ are distinct
points $x \neq y$ of $X$ then $\tau_A \otimes ev_x$ and $\tau_A \otimes ev_y$ are two distinct tracial states on
$A \otimes C(X)$. Here $ev_x$ is the functional on $C(X)$ given by $ev_x(f) \equiv f(x)$, $f \in C(X)$. 
\par
We record this as the following

\begin{thm}
The reduced $C^*$-algebra $C^*_r(BS(m,n))$ of the Baumslag-Solitar group $BS(m,n)$ is simple if and only if it
has a unique trace, if and only if $|n|, |m| \geq 2$ and $|n| \neq |m|$. 
\end{thm}
For more on $C^*$-simplicity of various groups see \cite{dlH05}. \\ 
\par

{\em Acknowledgements.} Most of the work on this note was done during my stay in Westf$\ddot{a}$lische
Wilhelms-Universit$\ddot{a}$t M$\ddot{u}$nster. I want to thank prof. J. Cuntz and prof. S. Echterhoff for their
hospitality. I want also to thank my advisor prof. Ken Dykema for some useful conversations I had with him.
Finally, I want to thank prof. Y. Ueda for some discussions on HNN-extensions.

\end{document}